\documentclass{article}

\usepackage{amsmath}
\usepackage{amssymb}
\usepackage{amsthm}
\usepackage[english]{babel}
\usepackage[T1]{fontenc}
\usepackage[utf8]{inputenc}
\usepackage{tikz}
\usepackage{xcolor}
\usepackage{comment}

\usepackage{latexsym,color,graphicx,enumerate}
\usepackage{hyperref}

\newtheorem{thm}{Theorem}[section]
\newtheorem{prop}[thm]{Proposition}
\newtheorem{lemma}[thm]{Lemma}
\newtheorem{rem}[thm]{Remark}

\makeatletter

\newcommand*{\plim}[1][]{%
	\if\relax\detokenize{#1}\relax
	\def\next{\qopname\relax m{lim}}%
	\else
	\def\next{\qopname\newmcodes@ m{#1-lim}}%
	\fi
	\next
}

\newcommand*{\psum}[1][]{%
	\DOTSB
	\if\relax\detokenize{#1}\relax\else
	\operatorname{#1-}\mkern-\thinmuskip
	\fi
	\sum@\slimits@
}

\makeatother

\newcommand{\R}{\mathbb{R}}             
\newcommand{\Z}{\mathbb{Z}}             
\newcommand{\C}{\mathbb{C}}             

\newcommand {\clr} {\color{red}}

\newcommand{\mn}{\vspace{0.2cm}\noindent}

\numberwithin{equation}{section}

\title{Inverse Spectral Analysis for the MIT Bag Model}

\author{Damien Gobin, Beno\^it Gr\'ebert, Bernard Helffer and Fran{\c{c}}ois Nicoleau}

\date{}

\begin{document}
	
	\maketitle

\begin{abstract}
We study an inverse spectral problem for radial Dirac operators with MIT bag
boundary conditions. We show that the spectra for effective angular momenta $\kappa=0$
and $\kappa=1$ uniquely determine both components of the potential in a
neighborhood of zero in $L^2(0,1)^2$. More precisely, the 
spectral map is a local real-analytic diffeomorphism at the zero potential.
\end{abstract}

\section{Introduction}

The MIT bag model was introduced in 1974 by Chodos et al.~\cite{Chodos1974}
to describe the confinement of quarks inside a hadron. In the static
spherical case, the hadron is represented by a ball of fixed radius.
Inside the ball, the quarks satisfy the free massless Dirac equation.
Confinement is imposed by a boundary condition on the spinor which
ensures that the normal component of the probability current vanishes
at the surface. Chodos et al.~\cite{ChodosBaryon1974} derived an equation
for the energy levels of a quark in the ball and computed the lowest
levels numerically. They then used these levels to predict the masses
of low-lying baryon states. The same confinement can also be obtained
by introducing a large mass outside the domain. Arrizabalaga
et al.~\cite{ALMR2019} showed that the eigenvalues converge to those
of the MIT bag operator as the exterior mass tends to infinity.

\medskip\noindent
By spherical symmetry, the Dirac equation reduces to a family of radial
systems. In this paper, we add a matrix potential with two unknown
components to these radial systems and keep the MIT bag boundary
condition. We aim to give an answer to the following question.
Do the spectra for two distinct effective angular momenta determine
both components of the potential?

\medskip\noindent
We studied the same question for radial Schr\"odinger operators
in~\cite{GGHN2025}. We proved that, for several pairs of angular
momenta, the two Dirichlet spectra determine the scalar potential
in a neighborhood of zero. In~\cite{GGHN2026}, we obtained a similar
local uniqueness result for singular radial AKNS operators with a
two-component potential and the boundary condition $Z_2(1)=0\,$.
Here we impose the MIT bag boundary condition, which in our convention
takes the form $Z_1(1)+Z_2(1)=0\,$.

\medskip\noindent
Following~\cite{GGHN2026}, we denote the \emph{effective angular momentum}
by $\kappa\,$. In the three-dimensional Dirac equation, $-\kappa$ is an
eigenvalue of the spin--orbit operator, so that $\kappa$ takes the values
$\pm 1,\pm 2,\ldots\,$. 
The regular case $\kappa=0$ arises naturally in the two-dimensional
MIT bag model on a disk with an Aharonov--Bohm flux at the origin
and the MIT bag boundary condition on the circle
\cite{BeneventanoEtAl2000}. When the flux equals half a flux quantum,
both $\kappa=0$ and $\kappa=1$ occur in the radial decomposition
(see~\cite[Appendix~A.2]{GGHN2026} and \cite{Th92}).

\medskip\noindent
On the interval $(0,1)$, the radial Dirac operator is given by
\begin{equation}\label{eq:AKNS}
    H_\kappa(V)\,Z = H_\kappa(p,q)\,Z
    =
    \begin{pmatrix}
        0 & -1\\[0.2em]
        1 & 0
    \end{pmatrix} Z'
    +
    \begin{pmatrix}
        0 & -\dfrac{\kappa}{x}\\[0.4em]
        -\dfrac{\kappa}{x} & 0
    \end{pmatrix} Z
    + V(x)\,Z\,,
    \qquad
    Z=(Z_1,Z_2)^{\mathsf T}\,,
\end{equation}
where $V$ is the matrix potential
\begin{equation}\label{eq:V}
    V(x)
    =
    \begin{pmatrix}
        -q(x) & p(x)\\[0.2em]
        p(x) & q(x)
    \end{pmatrix}\,,
    \qquad
    p,q \in L_\R^2(0,1)\,.
\end{equation}

\medskip\noindent
To describe the boundary conditions, we set
$u_\theta=(\sin\theta,\cos\theta)^{\mathsf T}$ for $\theta\in\mathbb R\,$.
When $\kappa=0$, both endpoints are regular, and we impose
\begin{equation}\label{eq:MITa}
    Z(0)\cdot u_{\theta_1}=0\,,
    \qquad
    Z(1)\cdot u_{\theta_2}=0\,.
\end{equation}
When $\kappa\neq0$, the origin is in the limit-point case, so we impose
a boundary condition only at $x=1$,
\begin{equation}\label{eq:MITb}
    Z(1)\cdot u_{\theta_2}=0\,.
\end{equation}
The MIT bag boundary conditions correspond to $\theta_2=\pm\pi/4\,$.
To see how the signs of $\kappa$ and $\theta_2$ are related, we use
the Pauli matrices
\begin{equation}\label{paulimatrices}
    \sigma_1=\begin{pmatrix}0&1\\[0.2em]1&0\end{pmatrix}\,,\quad
    \sigma_2=\begin{pmatrix}0&-i\\[0.2em]i&0\end{pmatrix}\,,\quad
    \sigma_3=\begin{pmatrix}1&0\\[0.2em]0&-1\end{pmatrix}\,.
\end{equation}
A direct computation gives the following identities for the differential
expressions:
\begin{align}
    \sigma_1 H_\kappa(p,q)\sigma_1 &= -H_{-\kappa}(-p,q)\,,
    \label{eq:sym-sigma1}\\
    \sigma_2 H_\kappa(p,q)\sigma_2 &= H_{-\kappa}(-p,-q)\,,
    \label{eq:sym-sigma2}\\
    \sigma_3 H_\kappa(p,q)\sigma_3 &= -H_\kappa(p,-q)\,.
    \label{eq:sym-sigma3}
\end{align}
At $x=1\,$, the two MIT bag boundary conditions are $Z_2(1)=-Z_1(1)$
and $Z_2(1)=Z_1(1)\,$. Applying $\sigma_1$ preserves each of these
relations, while applying $\sigma_2$ or $\sigma_3$ changes one into
the other. Thus, for $\kappa\neq0\,$, we may use $\sigma_1$ to change the sign
of $\kappa$ and $\sigma_3$ to change the sign of $\theta_2$. When $\kappa=0$, we choose $\theta_1=0\,$, that is, $Z_2(0)=0$.
If $Z$ satisfies this condition and $Z_1(1)-Z_2(1)=0 \,$, then
$\widetilde Z=\sigma_3 Z=(Z_1,-Z_2)^T$ satisfies
$\widetilde Z_2(0)=0$ and $\widetilde Z_1(1)+\widetilde Z_2(1)=0\,$.
Thus $\sigma_3$ changes $\theta_2=-\pi/4$ into $\theta_2=\pi/4$
while preserving the boundary condition at the origin. 

\medskip\noindent
In view of these symmetries, we may assume throughout the paper that
$ \kappa\geq 0$ and $\theta_2=\frac{\pi}{4}\,$.

\noindent
With this choice of boundary conditions, the domain of $H_0$ is
\begin{equation}\label{eq:Domaina}
    D(H_0)
    =
    \Bigl\{
    Z\in L^2(0,1)^2 \;:\;
    H_0 Z\in L^2(0,1)^2,\;
    Z_2(0)=0\,,\;
    Z_1(1)+Z_2(1)=0
    \Bigr\}\,.
\end{equation}
For $\kappa>0$, the domain is
\begin{equation}\label{eq:Domainb}
    D(H_\kappa)
    =
    \Bigl\{
    Z\in L^2(0,1)^2 \;:\;
    H_\kappa Z\in L^2(0,1)^2\,,\;
    Z_1(1)+Z_2(1)=0
    \Bigr\}\,.
\end{equation}

\medskip\noindent
For each $\kappa\geq0\,$, the operator $H_\kappa(p,q)$ is self-adjoint
and has purely discrete spectrum~\cite{Serier2006}. Its eigenvalues
are simple, and we denote them by
\begin{equation}\label{eq:EigenSeq}
    \lambda_{\kappa,n}(p,q)\,, \qquad n\in\mathbb Z\,,
\end{equation}
in increasing order. The labelling is fixed by the asymptotic formula
of~\cite[Theorem~3.1]{Serier2006},
\begin{equation}\label{eq:AsymptAKNS}
    \lambda_{\kappa,n}(p,q)
    =
    \left(n+\operatorname{sgn}(n)\frac{\kappa}{2}\right)\pi
    +\frac{\pi}{4}+r_{\kappa,n}\,,
    \qquad
    (r_{\kappa,n})_{n\in\mathbb Z}\in\ell^2(\mathbb Z)\,.
\end{equation}
\medskip\noindent
At the zero potential, the spectra for $\kappa=0$ and $\kappa=1$ can
be described as follows. For $\kappa=0\,$, the eigenvalues are
\[
    \lambda_{0,n}(0,0)=n\pi+\frac{\pi}{4}\,,
    \qquad n\in\mathbb Z\,.
\]
For $\kappa=1\,$, the eigenvalues $\lambda_{1,n}(0,0)$ are the nonzero
real roots of
\[
    J_{1/2}(\lambda)-J_{3/2}(\lambda)=0\,,
\]
or equivalently,
\[
    (\lambda-1)\sin\lambda+\lambda\cos\lambda=0\,.
\]
Both cases follow from the characteristic equation
\eqref{eq:Characteristic-MIT} derived below.

\medskip\noindent
For a fixed effective angular momentum $\kappa$, Serier~\cite{Serier2006}
showed that the eigenvalues together with suitable norming constants
determine the potential. The eigenvalues correspond to energy levels,
which can in principle be measured. The norming constants, however,
depend on the eigenfunctions and are not directly accessible from
measurements of these levels. This motivates the use of a second
spectrum. We keep the same MIT bag boundary condition and use the
eigenvalues for two distinct values of $\kappa\,$.

\medskip\noindent
To state our result, we introduce the renormalized eigenvalues
\[
    \widetilde{\lambda}_{\kappa,n}(p,q)
    =
    \lambda_{\kappa,n}(p,q)
    -
    \left(n+\operatorname{sgn}(n)\frac{\kappa}{2}\right)\pi
    -\frac{\pi}{4}\,,
    \qquad n\in\mathbb Z\,.
\]
By \eqref{eq:AsymptAKNS}, the sequence
$(\widetilde{\lambda}_{\kappa,n}(p,q))_{n\in\mathbb Z}$ belongs to
$\ell^2(\mathbb Z;\mathbb R)$. We set
\[
    \mathcal H=L^2((0,1);\mathbb R)\,,\qquad
    X=\mathcal H\times\mathcal H\,,\qquad
    Y=\ell^2(\mathbb Z;\mathbb R)\times\ell^2(\mathbb Z\,;\mathbb R)\,.
\]
For two distinct effective angular momenta $\kappa_1$ and $\kappa_2$,
we define the spectral map $\mathcal S_{\kappa_1,\kappa_2}:X\to Y$ by
\begin{equation}\label{eq:spectral-map-two-ell}
    \mathcal S_{\kappa_1,\kappa_2}(p,q)
    =
    \left(
        \bigl(\widetilde{\lambda}_{\kappa_1,n}(p,q)\bigr)_{n\in\mathbb Z},
        \bigl(\widetilde{\lambda}_{\kappa_2,n}(p,q)\bigr)_{n\in\mathbb Z}
    \right)\,.
\end{equation}

\medskip\noindent
Our main result is the following.

\begin{thm}\label{thm:main-MIT-01}
    There exist an open neighborhood $U$ of $(0,0)$ in $X$ and an open
    neighborhood $W$ of $\mathcal S_{0,1}(0,0)$ in $Y$ such that
    \[
        \mathcal S_{0,1}:U\longrightarrow W
    \]
    is a real-analytic diffeomorphism. In particular, the spectra for
    $\kappa=0$ and $\kappa=1$ uniquely determine both components of
    the potential $(p,q)$ in $U\,$.
\end{thm}

\medskip\noindent
The proof relies on the study of the Fr\'echet differential
\[
    D_{(p,q)}\mathcal S_{0,1}(0,0):X\longrightarrow Y\,.
\]
We first prove that this differential is injective. We then show that
it is Fredholm of index zero, so it is also surjective. Since the
spectral map is real-analytic by~\cite[Theorem~3.3]{Serier2006},
the result follows from the analytic inverse function theorem.

\medskip\noindent
For other pairs of effective angular momenta, the calculations needed
to study the differential are more involved. We expect that a new idea
will be needed to treat the general case.

\medskip\noindent
The paper is organized as follows. In Section~2, we describe the
eigenvalues and eigenfunctions at the zero potential. In Section~3,
we compute the differential of the spectral map at zero and formulate
the linearized inverse problem. In Section~4, we derive a Kneser--Sommerfeld formula for the MIT bag
boundary condition. Sections~5 and~6 are devoted to the proof of
Theorem~\ref{thm:main-MIT-01}.

\section{Eigenvalue analysis in the unperturbed case $V=0$\,}

In this section we analyze the spectral problem in the unperturbed case $V=0\,$,
which serves as the reference configuration for the perturbative and inverse
analysis developed later.

\medskip\noindent
When $V=0\,$, the AKNS operator reduces to the first-order matrix equation
\begin{equation}\label{ODEmatrix}
\begin{pmatrix}
0 & -1\vphantom{\dfrac{\kappa}{x}}\\[0.3em]
1 & 0
\end{pmatrix} Z'(x)
+
\begin{pmatrix}
0 & -\dfrac{\kappa}{x}\\[0.5em]
-\dfrac{\kappa}{x} & 0
\end{pmatrix} Z(x)
=\lambda Z(x)\,,
\qquad x\in(0,1)\,,
\end{equation}
with the boundary condition 
\begin{equation}\label{BC-V0}
	Z_2(0)=0\,,\; Z_1(1)+Z_2(1)=0
	\quad \text{for } \kappa = 0\,,
	\qquad
	Z_1(1)+Z_2(1)=0
	\quad \text{for } \kappa \neq 0\, .
\end{equation}

\subsection{The case $\lambda=0$}

The case $\lambda=0$ requires a specific discussion.  
Setting $\lambda=0$ in the unperturbed Dirac equation yields
\[
Z_1'=\frac{\kappa}{x}\,Z_1\,,
\qquad
Z_2'=-\frac{\kappa}{x}\,Z_2\, .
\]
Hence the general solutions are
\[
Z_1(x)=C_1\,x^{\kappa}\,,
\qquad
Z_2(x)=C_2\,x^{-\kappa}\,.
\]
For $\kappa=0\,$, the boundary condition $Z_2(0)=0$ forces $C_2=0$.
For $\kappa\geq 1\,$, the condition $Z_2\in L^2(0,1)$ near $x=0$
again forces $C_2=0\,$, so that
\[
Z(x)=\binom{C_1\,x^{\kappa}}{0}\,.
\]
The MIT bag boundary condition at $x=1$,
\[
Z_1(1)+Z_2(1)=0\,,
\]
then implies $C_1=0\,$. Hence there is no nontrivial solution for
$\lambda=0$, and therefore
\[
0\notin \operatorname{Spec}(H_\kappa(0))
\]
for every $\kappa\geq0\,$.

\subsection{The case $\lambda \not=0$}

We now consider the case $\lambda \neq 0$, for which the system \eqref{ODEmatrix} decouples into two scalar equations for $Z_1$ and $Z_2$:
	\begin{equation}\label{ODE-Z1Z2}
		\left\{
		\begin{aligned}
			&Z_1''(x)+\left(\lambda^2-\frac{\kappa(\kappa-1)}{x^2}\right)Z_1(x)=0\,,\\[0.4em]
			&Z_2''(x)+\left(\lambda^2-\frac{\kappa(\kappa+1)}{x^2}\right)Z_2(x)=0\,.
		\end{aligned}
		\right.
	\end{equation}
	We set $\nu=\kappa+\tfrac12$.
	After the standard substitution $Z_j(x)=\sqrt{x}\,u_j(\lambda x)\,$, the system reduces to Bessel equations of orders $\nu-1$ and $\nu$.
	Accordingly, a fundamental system of solutions is given by
	\begin{equation}\label{eq:Fundamental}
		Z^{(0)}(x,\lambda)
		=
		\sqrt{\frac{\pi \lambda x}{2}}\,
		\begin{pmatrix}
			J_{\nu-1}(\lambda x)\\[0.4em]
			-\,J_{\nu}(\lambda x)
		\end{pmatrix}\,,
		\qquad
		W^{(0)}(x,\lambda)
		=
		\sqrt{\frac{\pi \lambda x}{2}}\,
		\begin{pmatrix}
			-\,Y_{\nu-1}(\lambda x)\\[0.4em]
			Y_{\nu}(\lambda x)
		\end{pmatrix}\,.
	\end{equation}
where $J_\nu$ and $Y_\nu$ denote the Bessel functions of the first and
second kinds (see \cite{Lebedev72}).  We recall that for $z\in\mathbb{C}\setminus(-\infty,0]$ and any real or complex parameter
$\nu\,$, one has
\begin{equation}\label{BesselJ}
J_{\nu}(z)
=\Bigl(\frac{z}{2}\Bigr)^{\nu}
\sum_{k=0}^{\infty}
(-1)^{k}\,
\frac{\bigl(\tfrac{z^{2}}{4}\bigr)^{k}}{k!\,\Gamma(\nu+k+1)}\,.
\end{equation}
For $\nu\notin\mathbb{Z}\,$, the function $Y_\nu$ is given by \footnote{For an integer $n$, the function $Y_n(z)$ is defined by the limit
\(
Y_n(z)=\lim_{\nu\to n} Y_\nu(z)\,.
\)}

\begin{equation}\label{BesselY}
Y_{\nu}(z)
=
\frac{
   J_{\nu}(z)\,\cos(\nu\pi)\;-\; J_{-\nu}(z)
}{
   \sin(\nu\pi)
}\,.
\end{equation}
We emphasize that for half-integer orders \(\nu\in \tfrac12+\mathbb{N}\), the function
\(\sqrt{z}\,J_\nu(z)\) is entire, whereas \(\sqrt{z}\,Y_\nu(z)\) is
generally meromorphic because of its singularity at \(z=0\).

\vspace{0.2cm}\noindent
As $x\to 0$\,, these functions satisfy the classical asymptotics
:
\begin{equation}\label{singularite}
J_\nu(x)\sim \frac{1}{\Gamma(\nu+1)}
\left(\frac{x}{2}\right)^\nu\,,
\qquad
Y_\nu(x)\sim -\,\frac{\Gamma(\nu)}{\pi}
\left(\frac{2}{x}\right)^\nu\,,
\qquad \nu>0\,.
\end{equation}
Hence $J_\nu$ is regular at the origin, whereas $Y_\nu$ is singular.
For $\kappa\geq 1$, among the two fundamental solutions in
\eqref{eq:Fundamental}, only $Z^{(0)}(x,\lambda)$ is square--integrable
near $x=0$. In the case $\kappa=0\,$, both solutions are square--integrable
near $x=0\,$, but the boundary condition $Z_2(0)=0$ selects
$Z^{(0)}(x,\lambda)\,$.
For each eigenvalue $\lambda=\lambda_{\kappa,n}(0,0)$ of the unperturbed operator,
we define the associated eigenfunction by
\begin{equation}\label{eq:Fundamentalbis}
Z_{\kappa,n}^{(0)}(x)
=
c_{\kappa,n}
\begin{pmatrix}
Z_{1,\kappa,n}^{(0)}(x)\\[0.2em]
Z_{2,\kappa,n}^{(0)}(x)
\end{pmatrix}
=
c_{\kappa,n}\,
\sqrt{\lambda_{\kappa,n}(0,0) x}\,
\begin{pmatrix}
	J_{\nu-1}\!\bigl(\lambda_{\kappa,n}(0,0)\,x\bigr)\\[0.4em]
	-\,J_{\nu}\!\bigl(\lambda_{\kappa,n}(0,0)\,x\bigr)
\end{pmatrix}\,,
\end{equation}
where $c_{\kappa,n}>0$ is a normalization constant to be specified later.
Using the MIT bag boundary condition at $x=1\,$, we obtain that the
eigenvalues of $H_\kappa(0)$ are the zeros of the map
\[
\lambda\mapsto J_{\nu-1}(\lambda)-J_\nu(\lambda)\,,
\]
as will be detailed in the paragraph below.

\subsection{Symmetries}

In this paragraph we discuss the symmetry properties of the unperturbed
spectrum for $V=0\,$.

\medskip
\noindent
For the MIT bag boundary condition $\theta_2=\frac{\pi}{4}$, the
characteristic equation is
\begin{equation}\label{eq:Characteristic-MIT}
	J_{\nu-1}(\lambda)-J_\nu(\lambda)=0\,,
	\qquad
	\nu=\kappa+\frac12\,.
\end{equation}
In particular, the corresponding unperturbed spectrum is not symmetric with respect to zero. Indeed, let $Z$ be an eigenfunction for $\lambda$. By
\eqref{eq:sym-sigma3}, the solution selected at the origin for
$-\lambda$ is proportional to $\sigma_3 Z=(Z_1,-Z_2)^{\mathsf T}$.
Since $Z_1(1)=-Z_2(1)\neq 0$, this solution does not satisfy the
boundary condition at $x=1$. Hence $-\lambda$ is not an eigenvalue. This asymmetry can also be seen directly from \eqref{eq:AsymptAKNS},
which gives, for all sufficiently large $n$,
\[
    \lambda_{\kappa,-n-1}(0,0)
    < -\lambda_{\kappa,n}(0,0)
    < \lambda_{\kappa,-n}(0,0)\,.
\]

\subsection{Summary and notation}

For the unperturbed operator $H_\kappa(0,0)$, we write
\[
    \lambda_{\kappa,n}:=\lambda_{\kappa,n}(0,0)\,,
    \qquad n\in\mathbb Z\,.
\]
These eigenvalues are simple and are precisely the nonzero real roots
of \eqref{eq:Characteristic-MIT}. They are indexed in increasing order,
with the labelling fixed by \eqref{eq:AsymptAKNS} at $(p,q)=(0,0)\,$:
\[
    \lambda_{\kappa,n}
    =
    \left(n+\operatorname{sgn}(n)\frac{\kappa}{2}\right)\pi
    +\frac{\pi}{4}+r_{\kappa,n}\,,
    \qquad (r_{\kappa,n})_{n\in\mathbb Z}\in\ell^2(\mathbb Z)\,.
\]

\medskip\noindent
With $\nu=\kappa+\frac12$, the corresponding normalized eigenfunctions
are given by
\begin{equation}\label{eq:eigenfunctions-MIT}
    Z_{\kappa,n}^{(0)}(x)
    =
    c_{\kappa,n}\sqrt{\lambda_{\kappa,n}x}\,
    \begin{pmatrix}
        J_{\nu-1}(\lambda_{\kappa,n}x)\\[0.4em]
        -J_\nu(\lambda_{\kappa,n}x)
    \end{pmatrix}\,,
\end{equation}
where $c_{\kappa,n}>0$ is chosen so that
$\|Z_{\kappa,n}^{(0)}\|_{L^2(0,1)^2}=1\,$.

\vspace{0.5cm}\noindent
We conclude this section by giving explicit formulas for the eigenvalues
and eigenfunctions when $\kappa=0$.

\subsection{The case $\kappa=0$}\label{sec:k=0}

When $\kappa=0$, we have $\nu=\tfrac12$. The Bessel functions of orders $\pm\tfrac12$ are given by
\cite[Eq.~(10.16.1)]{DLMF}
\begin{equation}\label{eq:HalfInt}
    J_{-1/2}(z)=\sqrt{\frac{2}{\pi z}}\,\cos z\,,
    \qquad
    J_{1/2}(z)=\sqrt{\frac{2}{\pi z}}\,\sin z\,.
\end{equation}
Substituting these formulas into \eqref{eq:Fundamental}, we obtain
\[
    Z^{(0)}(x,\lambda)
    =
    \begin{pmatrix}
        \cos(\lambda x)\\[0.2em]
        -\sin(\lambda x)
    \end{pmatrix}\,,
\]
which satisfies $Z_2(0)=0$. The MIT bag boundary condition
$Z_1(1)+Z_2(1)=0$ then gives $\cos\lambda-\sin\lambda=0$. Hence
\begin{equation}\label{eq:l0-eigs}
    \lambda_{0,n}=n\pi+\frac{\pi}{4}\,,
    \qquad n\in\mathbb Z\,.
\end{equation}
The
corresponding eigenfunction already has unit norm in $L^2(0,1)^2\,$.
Thus $c_{0,n} = 1$ and
\begin{equation}\label{eq:l0-eigfct}
    Z_{0,n}^{(0)}(x)
    =
    \begin{pmatrix}
        \cos(\lambda_{0,n}x)\\[0.2em]
        -\sin(\lambda_{0,n}x)
    \end{pmatrix}\,.
\end{equation}

\section{Spectral map and the linearized problem at $V=0$}

In this section we introduce the spectral map and analyze its linearization at
the unperturbed configuration $V=0$. This linearized analysis provides the key
tool for the local inverse results established later.

\subsection{Differential of the eigenvalues and of the spectral map}

For each nonnegative integer $\kappa$ and each $n\in\mathbb Z$, the map
\[
    (p,q)\in X \longmapsto \lambda_{\kappa,n}(p,q)\in\mathbb R
\]
is real-analytic by \cite[Proposition~3.1]{Serier2006}.
Fix $(p,q)\in X$ and let
\[
    Z_{\kappa,n}(\, \cdot \, ;p,q)=(Z_1,Z_2)^{\mathsf T}
\]
be a real-valued normalized eigenfunction of $H_\kappa(p,q)$
associated with $\lambda_{\kappa,n}(p,q)\,$. For any
$v=(v_1,v_2)\in X$, the Fr\'echet differential is given by (see \cite[Proposition~3.1]{Serier2006})
\begin{equation}\label{eq:DiffGeneral}
    D_{(p,q)}\lambda_{\kappa,n}(p,q)\cdot v
    =
    \int_0^1
    \left[
        2Z_1(x)Z_2(x)v_1(x)
        +\bigl(Z_2(x)^2-Z_1(x)^2\bigr)v_2(x)
    \right]\,dx\,.
\end{equation}

\medskip\noindent
Since the term subtracted from $\lambda_{\kappa,n}(p,q)$ in the
renormalization is independent of $(p,q)$, the same formula holds
for $\widetilde{\lambda}_{\kappa,n}(p,q)\,$. Applying
\eqref{eq:DiffGeneral} for $\kappa=\kappa_1$ and $\kappa=\kappa_2$
therefore gives each component of the differential of the spectral
map $\mathcal S_{\kappa_1,\kappa_2}$ defined in
\eqref{eq:spectral-map-two-ell}.


\subsection{The linearized problem at $V=0$}

For two distinct effective angular momenta $\kappa_1$ and $\kappa_2\,$,
the linearized inverse problem is to determine whether
\[
    D_{(p,q)}\mathcal S_{\kappa_1,\kappa_2}(0,0)\cdot v=0
    \quad\Longrightarrow\quad v=0\,.
\]
Equivalently, if $v=(v_1,v_2)\in X$ satisfies
\[
    D_{(p,q)}\lambda_{\kappa_j,n}(0,0)\cdot v=0
    \quad \text{for all } n\in\mathbb Z \text{ and } j=1,2 \, ,
\]
must $v$ vanish?

\medskip\noindent
To write these conditions explicitly, fix $\kappa$ and $n$, and let
$Z_{\kappa,n}^{(0)}=(Z_1,Z_2)^{\mathsf T}$ be the normalized
eigenfunction of $H_\kappa(0,0)$ associated with $\lambda_{\kappa,n}$.
By \eqref{eq:DiffGeneral},
\begin{equation}\label{eq:LinEqGeneralBeta}
    D_{(p,q)}\lambda_{\kappa,n}(0,0)\cdot v
    =
    \int_0^1
    \left[
        2Z_1(x)Z_2(x)v_1(x)
        +\bigl(Z_2(x)^2-Z_1(x)^2\bigr)v_2(x)
    \right]\,dx\,.
\end{equation}
Substituting \eqref{eq:eigenfunctions-MIT} into this formula gives,
with $\nu=\kappa+\tfrac12\,$,
\begin{equation}\label{eq:DiffZeroSigned}
    \begin{aligned}
        D_{(p,q)}\lambda_{\kappa,n}(0,0)\cdot v
        = c_{\kappa,n}^{\,2}\int_0^1
        \lambda_{\kappa,n}x\,\Bigl[
        &-2J_{\nu-1}(\lambda_{\kappa,n}x)
            J_\nu(\lambda_{\kappa,n}x)\,v_1(x)\\
        &+\bigl(J_\nu(\lambda_{\kappa,n}x)^2
            -J_{\nu-1}(\lambda_{\kappa,n}x)^2\bigr)v_2(x)
        \Bigr]\,dx\,.
    \end{aligned}
\end{equation}
To show that the vanishing of these integrals for $\kappa=0$ and
$\kappa=1$ forces $v=0$, we shall use a Kneser--Sommerfeld formula,
which we derive in the next section.

\section{Kneser--Sommerfeld formula for the MIT bag model}

\subsection{Notation}

We fix a nonnegative integer $\kappa$, set $\nu=\kappa+\tfrac12$,
and write $\lambda_n=\lambda_{\kappa,n}(0,0)$. For negative real
arguments, we take boundary values from the upper half-plane and use
the principal branch of the square root. In particular,
\[
    J_\rho(-t)=e^{i\pi\rho}J_\rho(t)\,, \mbox{ for } t>0\,.
\]
For each fixed $n\in\mathbb Z$, the characteristic equation
\eqref{eq:Characteristic-MIT} allows us to write
\begin{equation}\label{eq:characteristic}
    J:=J_{\nu-1}(\lambda_n)=J_\nu(\lambda_n)\,.
\end{equation}

\medskip\noindent
Let $U_{\kappa,n}=c_{\kappa,n}^{-1}Z_{\kappa,n}^{(0)}$ denote the
unnormalized eigenfunction. By \eqref{eq:eigenfunctions-MIT},
\begin{equation}\label{eq:normalization-definition}
    c_{\kappa,n}^{-2}
    =\|U_{\kappa,n}\|_{L^2(0,1)^2}^2
    =\int_0^1 \lambda_n x
    \left[J_{\nu-1}(\lambda_n x)^2+J_\nu(\lambda_n x)^2\right]\,dx\,.
\end{equation}

\medskip\noindent
For the Kneser--Sommerfeld formula, we introduce the two functions
\begin{subequations}\label{eq:Phikappa}
\begin{align}
    \Phi_\kappa(z)
    &=
    \frac{\pi z}{2}
    \begin{pmatrix}
        -2J_{\nu-1}(z)J_\nu(z)\\[0.3em]
        J_\nu(z)^2-J_{\nu-1}(z)^2
    \end{pmatrix}\,,\\[0.5em]
    \Psi_\kappa(z)
    &=
    \frac{\pi z}{2}
    \begin{pmatrix}
        J_{\nu-1}(z)Y_\nu(z)+J_\nu(z)Y_{\nu-1}(z)\\[0.3em]
        J_{\nu-1}(z)Y_{\nu-1}(z)-J_\nu(z)Y_\nu(z)
    \end{pmatrix}\,.
\end{align}
\end{subequations}
With this notation, \eqref{eq:DiffZeroSigned} becomes
\begin{equation}\label{eq:differential-Phi}
    D_{(p,q)}\lambda_{\kappa,n}(0,0)\cdot v
    =
    \frac{2c_{\kappa,n}^2}{\pi}
    \int_0^1 \Phi_\kappa(\lambda_n x)\cdot v(x)\,dx\,,
    \qquad v\in X\,.
\end{equation}

\subsection{The normalization constants}

To compute $c_{\kappa,n}$, we use the Lommel formula
\cite[Eq.~(10.22.5)]{DLMF}. For $a>0$ and $\mu>-1$, it gives
\begin{equation}\label{eq:lommel-square}
    \int_0^1 xJ_\mu(ax)^2\,dx
    =
    \frac12\left[J_\mu(a)^2-J_{\mu-1}(a)J_{\mu+1}(a)\right]\,.
\end{equation}
We apply this formula with $\mu=\nu-1$ and $\mu=\nu$. Both orders
are greater than $-1$, including when $\kappa=0$. By our convention
for negative arguments, replacing $a$ by $-a$ multiplies both sides
of \eqref{eq:lommel-square} by $e^{2\pi i\mu}$. We may therefore
take $a=\lambda_n$ for every $n\in\mathbb Z\,$.

\medskip\noindent
Fix $n\in\mathbb Z$ and recall from \eqref{eq:characteristic} that
$J=J_{\nu-1}(\lambda_n)=J_\nu(\lambda_n)$. The recurrence relation
\cite[Eq.~(10.6.1)]{DLMF}
\begin{equation}\label{contiguite}
    J_{\rho-1}(z)+J_{\rho+1}(z)
    =\frac{2\rho}{z}J_\rho(z)\,, \qquad z\neq0\,,
\end{equation}
gives
\[
    J_{\nu-2}(\lambda_n)
    =\left(\frac{2(\nu-1)}{\lambda_n}-1\right)J\, , \qquad
    J_{\nu+1}(\lambda_n)
    =\left(\frac{2\nu}{\lambda_n}-1\right)J\,.
\]
Substituting into \eqref{eq:lommel-square}, we obtain
\[
    \begin{aligned}
        \int_0^1 xJ_{\nu-1}(\lambda_n x)^2\,dx
        &=\left(1-\frac{\nu-1}{\lambda_n}\right)J^2,\\
        \int_0^1 xJ_\nu(\lambda_n x)^2\,dx
        &=\left(1-\frac{\nu}{\lambda_n}\right)J^2\,.
    \end{aligned}
\]
Adding these identities and using $2\nu-1=2\kappa$, we find
\begin{equation}\label{eq:lommel-sum}
    \lambda_n\int_0^1 x
    \left[J_{\nu-1}(\lambda_n x)^2+J_\nu(\lambda_n x)^2\right]\,dx
    =2(\lambda_n-\kappa)J^2\,.
\end{equation}
By \eqref{eq:normalization-definition}, the left-hand side is
$\|U_{\kappa,n}\|_{L^2(0,1)^2}^2=c_{\kappa,n}^{-2}>0$\,.
In particular, $\lambda_n\neq\kappa$, and
\begin{equation}\label{eq:c-normalization}
    c_{\kappa,n}^2=\frac{1}{2(\lambda_n-\kappa)J^2}\,.
\end{equation}

\subsection{Definition of a meromorphic function}

Fix $z\in\mathbb C\setminus\{\lambda_n:n\in\mathbb Z\}\,$,
$v\in C_c^\infty((0,1);\mathbb C^2)$ and $x\in(0,1)$\,.
We introduce the meromorphic function
\begin{equation}\label{eq:r0}
	r_\kappa(w)
	=
	\frac{Y_{\nu-1}(w)-Y_\nu(w)}
	{J_{\nu-1}(w)-J_\nu(w)}
\end{equation}
and 
\begin{equation}\label{eq:F}
	\mathcal F(w;x,z)
	=
	-\frac{
		\left[\Psi_\kappa(wx)
		+r_\kappa(w)\Phi_\kappa(wx)\right]\cdot v(x)
	}{w-z}\,,
\end{equation}
where $\Phi_\kappa$ and $\Psi_\kappa$ were introduced in \eqref{eq:Phikappa}.

\medskip\noindent
For fixed $x\in(0,1)$, the function
\[
    w\longmapsto \Psi_\kappa(wx)+r_\kappa(w)\Phi_\kappa(wx)
\]
extends holomorphically across $w=0$. Indeed, both Bessel orders
$\nu-1=\kappa-\tfrac12$ and $\nu=\kappa+\tfrac12$ are half-integers.
The square-root factors cancel in the quotient defining $r_\kappa$
and in the expressions defining $\Phi_\kappa$ and $\Psi_\kappa\,$.
These functions are therefore meromorphic in $w$. As $w\to0$,
\[
    r_\kappa(w)=O(w^{-2\kappa})\,,\qquad
    \Phi_\kappa(wx)=O(w^{2\kappa})\,,\qquad
    \Psi_\kappa(wx)=O(1)\,.
\]
For $\kappa\geq1$, these estimates follow from \eqref{singularite}.
For $\kappa=0$, they follow directly from the elementary formulas
for Bessel functions of orders $\pm\tfrac12$.
It follows that $\Psi_\kappa(wx)+r_\kappa(w)\Phi_\kappa(wx)$
is bounded near $w=0$ and hence extends holomorphically across the origin.

\medskip\noindent
Since $z$ is outside the spectrum, $J_{\nu-1}(z)-J_\nu(z) \neq 0$ and $\mathcal F(\, \cdot\, ;x,z)$ has at most
a simple pole at $w=z\,$. By \eqref{eq:F}, its residue is
\begin{equation}\label{eq:residue-at-z}
    \operatorname*{Res}_{w=z}\mathcal F(w;x,z)
    =
    -\left[\Psi_\kappa(zx)+r_\kappa(z)\Phi_\kappa(zx)\right]\cdot v(x)\,.
\end{equation}
For $z=0$, the expression in brackets denotes the value at the origin
of the holomorphic extension obtained above.

\medskip\noindent
We next compute the residue of $r_\kappa$ at $\lambda_n$.
Recall from \eqref{eq:characteristic} that
$J=J_{\nu-1}(\lambda_n)=J_\nu(\lambda_n)$. The Wronskian identity
\begin{equation}\label{wronskian}
    J_\nu(t)Y_{\nu-1}(t)-J_{\nu-1}(t)Y_\nu(t)
    =\frac{2}{\pi t}
\end{equation}
gives, at $t=\lambda_n\,$,
\[
    J\left[Y_{\nu-1}(\lambda_n)-Y_\nu(\lambda_n)\right]
    =\frac{2}{\pi\lambda_n}\,.
\]
Thus $J\neq0$, and the numerator in \eqref{eq:r0} satisfies
\begin{equation}\label{eq:numerator}
    Y_{\nu-1}(\lambda_n)-Y_\nu(\lambda_n)
    =\frac{2}{\pi\lambda_n J}\,.
\end{equation}
To compute the derivative of $J_{\nu-1}(w)-J_\nu(w)$ at $w=\lambda_n\,$,
we use the recurrence relations \cite[Eq.~(5.3.5)]{Lebedev72}
\[
    J_\alpha'(t)=\frac{\alpha}{t}J_\alpha(t)-J_{\alpha+1}(t)\,,
    \qquad
    J_\alpha'(t)=J_{\alpha-1}(t)-\frac{\alpha}{t}J_\alpha(t)\,.
\]
At $t=\lambda_n$, these give
\[
    J_{\nu-1}'(\lambda_n)
    =\left(\frac{\nu-1}{\lambda_n}-1\right)J\,,
    \qquad
    J_\nu'(\lambda_n)
    =\left(1-\frac{\nu}{\lambda_n}\right)J\,.
\]
Since $2\nu-1=2\kappa\,$, we obtain
\begin{equation}\label{eq:denominator}
    J_{\nu-1}'(\lambda_n)-J_\nu'(\lambda_n)
    =\frac{2(\kappa-\lambda_n)}{\lambda_n}J\,.
\end{equation}
This derivative is nonzero since $\lambda_n\neq\kappa$ by
\eqref{eq:lommel-sum}. Thus the denominator in \eqref{eq:r0} has
a simple zero at $\lambda_n$, and the numerator is nonzero by
\eqref{eq:numerator}. It follows that $r_\kappa$ has a simple pole
at $\lambda_n$, with residue
\begin{align}
    \operatorname*{Res}_{w=\lambda_n}r_\kappa(w)
    &=
    \frac{Y_{\nu-1}(\lambda_n)-Y_\nu(\lambda_n)}
         {J_{\nu-1}'(\lambda_n)-J_\nu'(\lambda_n)}
    \label{eq:residue-r}\\
    &=-\frac{1}{\pi(\lambda_n-\kappa)J^2}
    \label{eq:residue-explicit}\\
    &=-\frac{2c_{\kappa,n}^2}{\pi}\,,
    \label{eq:r-residue}
\end{align}
where the last equality follows from \eqref{eq:c-normalization}. Finally, $\Phi_\kappa(wx)$ and $\Psi_\kappa(wx)$ are holomorphic
at $w=\lambda_n$. Since $\lambda_n\neq z\,$, \eqref{eq:F} gives
\begin{align}
    \operatorname*{Res}_{w=\lambda_n}\mathcal F(w;x,z)
    &=
    -\frac{\Phi_\kappa(\lambda_n x)\cdot v(x)}{\lambda_n-z}
    \operatorname*{Res}_{w=\lambda_n}r_\kappa(w)
    \notag\\
    &=
    \frac{2c_{\kappa,n}^2}{\pi(\lambda_n-z)}
    \Phi_\kappa(\lambda_n x)\cdot v(x)\,.
    \label{eq:residue-F}
\end{align}

\subsection{Contour-integral representation}

We keep the notation of the preceding subsection and set
$K=\operatorname{supp}v$. For each integer $N\geq1$, let
\[
    a_N=\frac{\lambda_N+\lambda_{N+1}}{2}\,,
    \qquad
    b_N=\frac{\lambda_{-N-1}+\lambda_{-N}}{2}\,,
    \qquad
    T_N=\pi N.
\]
Since $a_N\to+\infty\,$, $b_N\to-\infty\,$, and $T_N\to+\infty\,$, for each $z\in\C$\,,
there exists $N_0=N_0(z)$ such that, for every $N\geq N_0$,
\[
    b_N<\operatorname{Re}z<a_N\,,
    \qquad
    |\operatorname{Im}z|<T_N\,.
\]
Let $\Gamma_N$ be the positively oriented boundary
of the rectangle
\[
R_N=
\left\{
w\in\mathbb C:
b_N<\operatorname{Re}w<a_N\,,\quad
|\operatorname{Im}w|<T_N
\right\}\,.
\]
The eigenvalue asymptotics recalled in \eqref{eq:AsymptAKNS} give as $N$ tends to $+\infty$
\[
a_N-\lambda_N
=\lambda_{N+1}-a_N
=\frac{\pi}{2}+o(1)
\]
and
\[
\lambda_{-N}-b_N
=b_N-\lambda_{-N-1}
=\frac{\pi}{2}+o(1)\,.
\]
Thus, for all sufficiently large $N\,$, the vertical sides of $\Gamma_N$ remain
uniformly separated from the spectrum,
 the rectangle
contains $z$ and precisely the eigenvalues
$\lambda_n$ with $-N\le n\le N\,$ and 
its boundary contains no poles of $\mathcal F(\, \cdot \,;x ,z)\,$.

\begin{figure}[ht]
	\centering
	\begin{tikzpicture}[scale=0.95,>=stealth]
		\draw[->,gray] (-5,0) -- (5,0)
		node[right] {$\operatorname{Re}w$};
		\draw[->,gray] (0,-2.8) -- (0,2.9)
		node[above] {$\operatorname{Im}w$};
		
		\draw[thick,blue]
		(-3.4,-2) rectangle (3.4,2);
		
		\draw[->,thick,blue]
		(-1.8,-2) -- (-0.8,-2);
		\draw[->,thick,blue]
		(3.4,-1.3) -- (3.4,-0.6);
		\draw[->,thick,blue]
		(1.8,2) -- (0.8,2);
		\draw[->,thick,blue]
		(-3.4,1.3) -- (-3.4,0.6);
		
		\node[blue,above right] at (3.4,2)
		{$\Gamma_N$};
		\node[above left] at (-3.4,0) {$b_N$};
		\node[above right] at (3.4,0) {$a_N$};
		\node[above left] at (0,2) {$iT_N$};
		\node[below left] at (0,-2) {$-iT_N$};
		
		\fill (-2.5,0) circle (2pt);
		\node[below] at (-2.5,0) {$\lambda_{-N}$};
		\fill (2.5,0) circle (2pt);
		\node[below] at (2.5,0) {$\lambda_N$};
		\node[above] at (-1.3,0) {$\cdots$};
		\node[above] at (1.3,0) {$\cdots$};
		
		\fill (-4.3,0) circle (2pt);
		\node[below] at (-4.3,0)
		{$\lambda_{-N-1}$};
		\fill (4.3,0) circle (2pt);
		\node[below] at (4.3,0)
		{$\lambda_{N+1}$};
		
		\fill (1,0.9) circle (2pt);
		\node[above right] at (1,0.9) {$z$};
	\end{tikzpicture}
	\caption{The positively oriented contour $\Gamma_N$.
	}
	\label{fig:contour}
\end{figure}
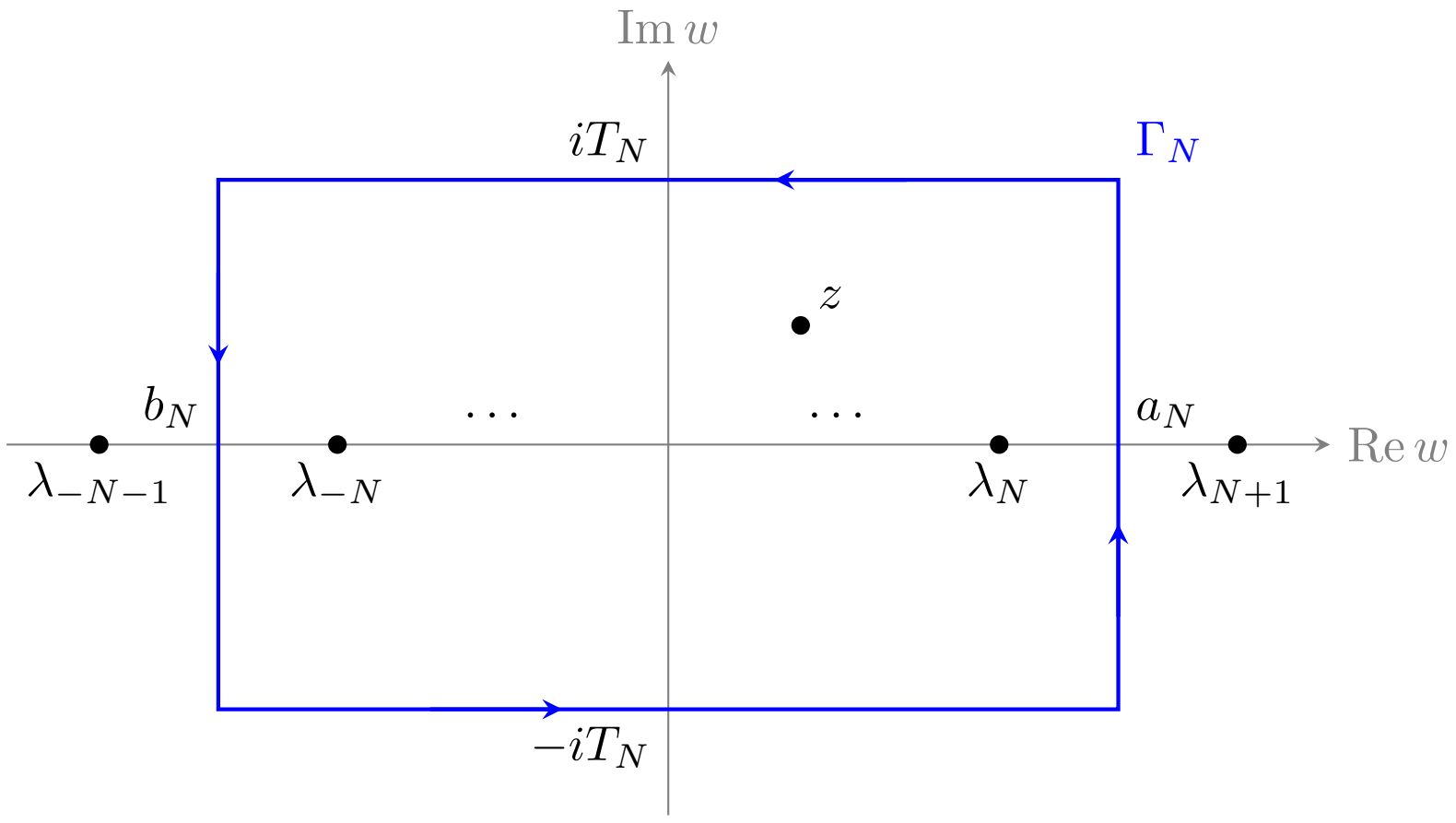

\medskip\noindent
It follows that for each fixed $x\in(0,1)\,$ and $z$ avoiding the spectrum, the residue theorem gives
\[
\frac{1}{2\pi i}
\int_{\Gamma_N}\mathcal F(w;x,z)\,dw
=
\operatorname*{Res}_{w=z}\mathcal F(w;x,z)
+
\sum_{n=-N}^{N}
\operatorname*{Res}_{w=\lambda_n}\mathcal F(w;x,z)\,.
\]
Using \eqref{eq:residue-at-z} and
\eqref{eq:residue-F}, we obtain the finite identity
\begin{align}
	&\left[
	\Psi_\kappa(zx)+r_\kappa(z)\Phi_\kappa(zx)
	\right]\cdot v(x)
	\notag\\
	&\quad=
	\sum_{n=-N}^{N}
	\frac{2c_{\kappa,n}^2}{\pi(\lambda_n-z)}
	\Phi_\kappa(\lambda_nx)\cdot v(x)
	-\frac{1}{2\pi i}
	\int_{\Gamma_N}\mathcal F(w;x,z)\,dw\,.
	\label{eq:finite-KS}
\end{align}

\medskip\noindent
With $z$ fixed, we show below that
\begin{equation}\label{contourvanish}
	\lim_{N\to\infty}
	\sup_{x\in K}
	\left|
	\int_{\Gamma_N}\mathcal F(w;x,z)\,dw
	\right|
	=0\,.
\end{equation}
Passing to the limit in \eqref{eq:finite-KS} yields
\begin{equation}\label{eq:pointwise-KS}
	\begin{aligned}
		&\left[
		\Psi_\kappa(zx)+r_\kappa(z)\Phi_\kappa(zx)
		\right]\cdot v(x)\\
		&\qquad=
		\lim_{N\to\infty}
		\sum_{n=-N}^{N}
		\frac{2c_{\kappa,n}^2}{\pi(\lambda_n-z)}
		\Phi_\kappa(\lambda_nx)\cdot v(x)\,.
	\end{aligned}
\end{equation}
The convergence is uniform for $x\in K$, and all terms vanish outside $K$.
We may therefore integrate
\eqref{eq:pointwise-KS} term by term over $(0,1)$.
Using \eqref{eq:differential-Phi}, we obtain
\begin{equation}\label{eq:kneser-sommerfeld}
	\int_0^1
	\left[\Psi_\kappa(zx)+r_\kappa(z)\Phi_\kappa(zx)\right]
	\cdot v(x)\,dx
	=
	\lim_{N\to\infty}\sum_{n=-N}^{N}
	\frac{D_{(p,q)}\lambda_{\kappa,n}(0,0)\cdot v}
	{\lambda_n-z}\,.
\end{equation}
Equation~\eqref{eq:kneser-sommerfeld} is the Kneser--Sommerfeld identity
for the MIT bag boundary condition. It holds for each fixed
$z\in\mathbb C\setminus\{\lambda_n:n\in\mathbb Z\}$, with the series
taken in the sense of symmetric partial sums. At $z=0$, the expression
in brackets is understood by holomorphic continuation. To complete
the proof, it remains to establish \eqref{contourvanish}.

\subsection{Estimates on the contour integral}\label{subsec:supcompact}

We prove \eqref{contourvanish} for fixed $\kappa$ and $z$.
Choose $\varepsilon\in(0,1/2)$ such that
\[
    K=\operatorname{supp}v\subset[\varepsilon,1-\varepsilon].
\]

\medskip\noindent
We shall prove that there exist a constant $C>0$ and an integer $N_1$,
both depending only on $\kappa$ and $\varepsilon$, such that
\begin{equation}\label{eq:KS-kernel-estimate}
    \left|\Psi_\kappa(wx)+r_\kappa(w)\Phi_\kappa(wx)\right|
    \leq C\left(
        e^{-2\varepsilon|\operatorname{Im}w|}
        +\frac{1}{|w|}
    \right)\,,
\end{equation}
for every $N\geq N_1\,$, every $w\in\Gamma_N\,$, and every
$x\in[\varepsilon,1-\varepsilon]\,$.
The proof uses the cancellation of the leading terms in the mixed
Hankel products of orders $\nu-1$ and $\nu\,$.

\medskip\noindent
\textbf{Normalized Hankel functions.}
For $\alpha\in\{\nu-1,\nu\}\,$, set
\[
h_\alpha^+(t)
=
\sqrt{\frac{\pi t}{2}}\,H_\alpha^{(1)}(t)\,,\, h_\alpha^-(t)
=
\sqrt{\frac{\pi t}{2}}\,H_\alpha^{(2)}(t)\,,
\]
where
\begin{equation} \label{eq:ident}
   Y_\alpha-iJ_\alpha=-iH_\alpha^{(1)}\,.
\qquad
Y_\alpha+iJ_\alpha=iH_\alpha^{(2)}\,, 
\end{equation}
Since the orders are half-integers, these normalized
functions extend meromorphically to $\mathbb C$,
with a possible pole only at the origin. More precisely, the half-integer Hankel formulas
give, for $t\in \mathbb C \setminus \{0\}\,$,
\begin{equation}\label{eq:KS-finite-Hankel}
	h_\alpha^\pm(t)
	=
	e^{\pm i(t-\alpha\pi/2-\pi/4)}
	P_\alpha^\pm(t^{-1})\,,
\end{equation}
where $P_\alpha^\pm$ are polynomials satisfying
$P_\alpha^\pm(0)=1\,$.
For instance, for $\alpha=-1/2$, both polynomials are identically
equal to one.
See \cite[Sections 10.17(i) and 10.49(i)]{DLMF}. Since $P_\alpha^\pm$ is a polynomial with
$P_\alpha^\pm(0)=1\,$, we have
\[
P_\alpha^\pm(t^{-1})=1+O(|t|^{-1})\,,
\]
uniformly for all complex $t$ with $|t|\ge1\,$.
Thus the resulting Hankel estimates apply
on all four sides of $\Gamma_N\,$, including
near the negative real axis.

\medskip\noindent
Set
\[
\theta(t)=t-\frac{\kappa\pi}{2}\,,
\qquad
\mathbf h^\pm(t)
=
\begin{pmatrix}
	h_{\nu-1}^\pm(t)\\
	h_\nu^\pm(t)
\end{pmatrix}\,.
\]
Since $\nu=\kappa+1/2\,$, we have, uniformly as
$|t|\to\infty\,$,
\begin{equation}\label{eq:KS-Hankel-vectors}
	\mathbf h^\pm(t)
	=
	e^{\pm i\theta(t)}
	\left[
	\begin{pmatrix}
		1\\
		\mp i
	\end{pmatrix}
	+O(|t|^{-1})
	\right]\,.
\end{equation}

\medskip\noindent
\textbf{A lower bound for the characteristic function.}
Introduce the normalized characteristic function
\[
d_\kappa(w)
=
\sqrt{\frac{\pi w}{2}}
\left[J_{\nu-1}(w)-J_\nu(w)\right]\,.
\]
By \eqref{eq:KS-Hankel-vectors} and the identity
$J_\alpha=(H_\alpha^{(1)}+H_\alpha^{(2)})/2$\,,
\begin{equation}\label{eq:KS-characteristic-asymptotics}
	d_\kappa(w)
	=
	\sqrt{2}\cos\left(
	w-\frac{\kappa\pi}{2}+\frac{\pi}{4}
	\right)
	+
	O\left(
	\frac{e^{|\operatorname{Im}w|}}{|w|}
	\right)\,.
\end{equation}
By the eigenvalue asymptotics \eqref{eq:AsymptAKNS} and the definitions
of $a_N$ and $b_N$\,, we have
\[
\operatorname{dist}\left(
\xi_N-\frac{\kappa\pi}{2}+\frac{\pi}{4},
\pi\mathbb Z
\right)=o(1)\,,
\qquad \xi_N\in\{a_N,b_N\}\,,
\mbox{ as } N\to\infty\,.
\]
Since, for real $u,y\,$,
\[
|\cos(u+iy)|^2=\cos^2u+\sinh^2y\,,
\]
it follows that, on the vertical sides, there exists a positive $c$ such that
\[
\left|
\cos\left(
w-\frac{\kappa\pi}{2}+\frac{\pi}{4}
\right)
\right|
\ge c \, e^{|\operatorname{Im}w|}\,.
\]
The same bound holds on the horizontal sides,
where $|\operatorname{Im}w|=T_N\to\infty\,$,
because $|\cos(u+iy)|\ge|\sinh y|\,$. Hence there
exist a constant $c>0$ and an integer $N_2\geq1$, both depending
only on $\kappa$, such that, for every $N\geq N_2$ and every
$w\in\Gamma_N$,
\begin{equation}\label{eq:KS-denominator-bound}
    |d_\kappa(w)|\geq c\,e^{|\operatorname{Im}w|}\,.
\end{equation}
Moreover, the definition \eqref{eq:r0} and the identities \eqref{eq:ident} 
yield
\[
r_\kappa(w)-i
=
-i\,
\frac{h_{\nu-1}^+(w)-h_\nu^+(w)}
{d_\kappa(w)}\,,
\]
and
\[
r_\kappa(w)+i
=
i\,
\frac{h_{\nu-1}^-(w)-h_\nu^-(w)}
{d_\kappa(w)}\,.
\]
Combining these identities with
\eqref{eq:KS-Hankel-vectors} and
\eqref{eq:KS-denominator-bound}, we obtain the existence of a positive $C$ such that:
\begin{equation}\label{eq:KS-r-exponential}
	\begin{aligned}
		|r_\kappa(w)-i|
		&\le C \, e^{-2y}\,,
		&& y=\operatorname{Im}w\ge0\,,\\
		|r_\kappa(w)+i|
		&\le C\, e^{-2|y|}\,,
		&& y=\operatorname{Im}w\le0\,,
	\end{aligned}
	\qquad w\in\Gamma_N\,.
\end{equation}
In particular, these estimates imply
\begin{equation}\label{eq:KS-r-bound}
	|r_\kappa(w)|\le C\,,
	\mbox{ for } w\in\Gamma_N\,.
\end{equation}

\medskip\noindent
\textbf{Cancellation in the vector kernel.}
We introduce the symmetric bilinear map
\[
\mathbb C^2 \ni (u,\widetilde u) \mapsto 
\mathcal B(u,\widetilde u)
=
\begin{pmatrix}
	u_1\widetilde u_2+u_2\widetilde u_1\\
	u_1\widetilde u_1-u_2\widetilde u_2
\end{pmatrix}\,,
\]
and the normalized Bessel vectors
\[
\mathbf j(t)
=
\sqrt{\frac{\pi t}{2}}
\begin{pmatrix}
	J_{\nu-1}(t)\\
	J_\nu(t)
\end{pmatrix}\,,
\qquad
\mathbf y(t)
=
\sqrt{\frac{\pi t}{2}}
\begin{pmatrix}
	Y_{\nu-1}(t)\\
	Y_\nu(t)
\end{pmatrix}\,.
\]
The definitions of $\Phi_\kappa$ and $\Psi_\kappa$ (see \eqref{eq:Phikappa})
give
\[
\Phi_\kappa(t)=-\mathcal B(\mathbf j(t),\mathbf j(t))\,,
\qquad
\Psi_\kappa(t)=\mathcal B(\mathbf j(t),\mathbf y(t))\,.
\]
For brevity, we write
\[
\mathcal Q_\kappa(w,x)
=
\Psi_\kappa(wx)+r_\kappa(w)\Phi_\kappa(wx)\,.
\]
Since
\[
\mathbf j=\frac{\mathbf h^++\mathbf h^-}{2}\,,
\qquad
\mathbf y=\frac{\mathbf h^+-\mathbf h^-}{2i}\,,
\]
the bilinearity and symmetry of $\mathcal B$ gives, with  $t=wx\,,$ the exact identity
\begin{equation}\label{eq:KS-kernel-Hankel}
	\begin{aligned}
		4\mathcal Q_\kappa(w,x)
		={}&(-i-r_\kappa(w))
		\mathcal B(\mathbf h^+(t),\mathbf h^+(t))\\
		&+(i-r_\kappa(w))
		\mathcal B(\mathbf h^-(t),\mathbf h^-(t))\\
		&-2r_\kappa(w)
		\mathcal B(\mathbf h^+(t),\mathbf h^-(t))\,.
	\end{aligned}
\end{equation}
We now use the crucial cancellation
\[
\mathcal B\left(
\begin{pmatrix}1\\-i\end{pmatrix},
\begin{pmatrix}1\\i\end{pmatrix}
\right)= \begin{pmatrix}0 \\0\end{pmatrix}\,.
\]
Thus \eqref{eq:KS-Hankel-vectors} implies 
\begin{equation}\label{eq:KS-mixed-cancellation}
	\left|
	\mathcal B(\mathbf h^+(t),\mathbf h^-(t))
	\right|
	\le \frac{C}{|t|}\,,
	\qquad |t|\ge1
\end{equation}
since the exponential factors cancel exactly in this
mixed product. The pure products satisfy
\begin{equation}\label{eq:KS-pure-products}
	\begin{aligned}
		\left|
		\mathcal B(\mathbf h^+(t),\mathbf h^+(t))
		\right|
		&\le C e^{-2\operatorname{Im}t}\,,\\
		\left|
		\mathcal B(\mathbf h^-(t),\mathbf h^-(t))
		\right|
		&\le C e^{2\operatorname{Im}t}\,,
	\end{aligned}
	\qquad |t|\ge1\,.
\end{equation}

\medskip\noindent
Now, suppose first that $y=\operatorname{Im}w\ge0$.
Using \eqref{eq:KS-r-exponential}, \eqref{eq:KS-r-bound},
\eqref{eq:KS-mixed-cancellation} and
\eqref{eq:KS-pure-products} in
\eqref{eq:KS-kernel-Hankel}, we find
\[
|\mathcal Q_\kappa(w,x)|
\le C\left(
e^{-2xy}
+e^{-2(1-x)y}
+\frac{1}{|w|x}
\right)\,.
\]
For $y\le0\,$, the roles of $\mathbf h^+$ and
$\mathbf h^-$ are interchanged, and the same
estimate holds with $y$ replaced by $|y|\,$.
Since $x\in[\varepsilon,1-\varepsilon]\,$,
we obtain \eqref{eq:KS-kernel-estimate}.
The condition $|wx|\ge1$ holds uniformly on
$\Gamma_N\times K$ for all sufficiently large $N$.

\medskip\noindent
\textbf{Integration over the four sides.}
The eigenvalue asymptotics imply
\[
a_N=\pi N+O(1)\,,\qquad
b_N=-\pi N+O(1)\,,\qquad
T_N=\pi N\,.
\]
Consequently, there exist constants $c,C>0$ and an integer
$N_0\geq1$, depending only on $\kappa$, such that
\[
    cN\leq |w|\leq CN
    \qquad\text{for every } N\geq N_0
    \text{ and every } w\in\Gamma_N\,.
\]
Since $z$ is fixed, for every $N\geq\max\{N_0,2|z|/c\}$ and every
$w\in\Gamma_N$ we have
\[
    |w-z|\geq |w|-|z|\geq \frac{cN}{2}\,.
\]
For $N\geq\max\{N_0,N_1,2|z|/c\}$, it follows from
\eqref{eq:F} and \eqref{eq:KS-kernel-estimate} that
\begin{equation}\label{eq:KS-F-contour-bound}
    \sup_{x\in K}|\mathcal F(w;x,z)|
    \leq \frac{C}{N}
    \left(e^{-2\varepsilon|\operatorname{Im}w|}+\frac1N\right)\,,
    \qquad w\in\Gamma_N\,.
\end{equation}
Let $\Gamma_N^{\mathrm h}$ and
$\Gamma_N^{\mathrm v}$ denote the unions of the
horizontal and vertical sides, respectively,
with their induced orientations.
On the horizontal sides,
$|\operatorname{Im}w|=\pi N\,$, and their total
length is $O(N)$. Hence
\[
\sup_{x\in K}
\left|
\int_{\Gamma_N^{\mathrm h}}\mathcal F(w;x,z)\,dw
\right|
\le C\left(e^{-2\varepsilon\pi N}+\frac1N\right)\,.
\]
On each vertical side, write $w=a_N+iy$
or $w=b_N+iy\,$, with $|y|\le T_N$.
Using \eqref{eq:KS-F-contour-bound}, we obtain
\[
\begin{aligned}
	\sup_{x\in K}
	\left|
	\int_{\Gamma_N^{\mathrm v}}\mathcal F(w;x,z)\,dw
	\right|
	&\le
	\frac{C}{N}
	\int_{-T_N}^{T_N}
	\left(e^{-2\varepsilon|y|}+\frac1N\right)\,dy\\
	&\le \frac{\tilde{C}}{N}\,.
\end{aligned}
\]
Combining the two estimates yields
\[
\sup_{x\in K}
\left|
\int_{\Gamma_N}\mathcal F(w;x,z)\,dw
\right|
\le
C \, \left(e^{-2\varepsilon\pi N}+\frac1N\right)
\longrightarrow0\,.
\]
This proves \eqref{contourvanish} and hence the Kneser--Sommerfeld
identity \eqref{eq:kneser-sommerfeld} for every
$v\in C_c^\infty((0,1);\mathbb C^2)$ and every fixed
$z\in\mathbb C\setminus\{\lambda_n:n\in\mathbb Z\}$.


\subsection{Extension to $L^2$ perturbations}

We extend \eqref{eq:pointwise-KS}-\eqref{eq:kneser-sommerfeld} to
$v\in L^2((0,1);\mathbb R^2)$ by proving
the convergence in $L^2$ of the partial-sum kernels
\[
S_N(z,x)
=
\sum_{n=-N}^{N}
\frac{2c_{\kappa,n}^2}{\pi(\lambda_n-z)}
\Phi_\kappa(\lambda_nx)\,.
\]
Throughout this subsection, $\kappa$ and
$z\in\mathbb C\setminus\{\lambda_n:n\in\mathbb Z\}$
are fixed.

\medskip\noindent
\textbf{Comparison with Fourier series.}
For real $t$ with $|t|\ge 1\,$, the half-integer
Bessel formula \cite[Eq.~(10.17.3)]{DLMF} gives
\[
\Phi_\kappa(t)
=
(-1)^\kappa
\begin{pmatrix}
	-\sin(2t)\\
	-\cos(2t)
\end{pmatrix}
+O(|t|^{-1})\,.
\]
Since $\Phi_\kappa$ is continuous at the origin,
we can write
\[
\Phi_\kappa(\lambda_nx)
=
(-1)^\kappa
\begin{pmatrix}
	-\sin(2\lambda_nx)\\
	-\cos(2\lambda_nx)
\end{pmatrix}
+R_n(x)\,,
\]
where
\[
|R_n(x)|
\le C\min\left\{1,\frac{1}{|\lambda_n|x}\right\}\,.
\]
Splitting the integral at $x=|\lambda_n|^{-1}\,$,
we obtain, for $|\lambda_n|\ge1\,$,
\[ 
	\|R_n\|_{L^2}^2
	\le
	C\int_0^{|\lambda_n|^{-1}}dx
	+\frac{C}{|\lambda_n|^2}
	\int_{|\lambda_n|^{-1}}^1\frac{dx}{x^2}\\
	 \le \frac{\widehat C}{|\lambda_n|}\,.
\]
Thus
\begin{equation}\label{eq:L2-KS-remainder}
	\|R_n\|_{L^2}=O(|n|^{-1/2})\,,
	\mbox{ as }  |n|\to\infty\,.
\end{equation}
The eigenvalue asymptotics \eqref{eq:AsymptAKNS} give
\[
\lambda_n=\pi n+\beta_\pm+\delta_n\,,
\qquad n\to\pm\infty\,,
\qquad (\delta_n)_{n\in\mathbb Z}\in\ell^2(\mathbb Z)\,,
\]
with real constants $\beta_\pm\,$.
Moreover, the Bessel asymptotics, the characteristic
equation \eqref{eq:characteristic} and \eqref{eq:c-normalization} yield
\[
c_{\kappa,n}^2=\frac{\pi}{2}+O(|n|^{-1})\,.
\]
Since $(\delta_n)$ is bounded, we have
$\lambda_n=\pi n+O(1)\,$, and hence 
\[
\frac{2c_{\kappa,n}^2}{\pi(\lambda_n-z)}
=
\frac{1}{\pi n}+O(|n|^{-2})\,.
\]
On the other hand,
\[
\sup_{0\le x\le1}
\left|
e^{2i\lambda_nx}
-e^{2i(\pi n+\beta_\pm)x}
\right|
\le 2|\delta_n|\,.
\]
Combining these estimates with
\eqref{eq:L2-KS-remainder}, we obtain
\begin{equation}\label{eq:L2-KS-Fourier-comparison}
	\frac{2c_{\kappa,n}^2}{\pi(\lambda_n-z)}
	\Phi_\kappa(\lambda_nx)
	=
	\frac{(-1)^\kappa}{\pi n}
	\begin{pmatrix}
		-\sin(2\pi nx+2\beta_\pm x)\\
		-\cos(2\pi nx+2\beta_\pm x)
	\end{pmatrix}
	+E_n(z,x)\,,
\end{equation}
where, for $|n|$ sufficiently large,
\[
\|E_n(z,\cdot)\|_{L^2}
\le C\left(
|n|^{-3/2}+\frac{|\delta_n|}{|n|}
\right)\,.
\]
Since $(\delta_n)\in\ell^2\,$, Cauchy--Schwarz inequality 
shows that the series of remainders converges
absolutely in $L^2\,$.

\medskip\noindent
Let $v_n$ denote the trigonometric vector in the leading term.
For $n$ and $m$ of the same sign, the phase shift is common
to both vectors and cancels in the inner product:
\[
\langle v_n,v_m\rangle_{L^2}
= \int_0^1 \cos(2\pi(n-m)x)\,dx
= \delta_{nm}\,.
\]
Thus the vectors form an orthonormal family on each tail.
Since $\sum_{n\ne 0} 1/n^2 < \infty\,$, the series of leading
terms converges in $L^2$ on each tail. Adding the two limits
gives convergence of the full leading-term series in $L^2\,$. Consequently, there exists
$S(z,\cdot)\in L^2((0,1);\mathbb C^2)$ such that
\begin{equation}\label{eq:L2-KS-kernel-convergence}
	S_N(z,\cdot)\longrightarrow S(z,\cdot)
	\quad\text{in }L^2((0,1);\mathbb C^2)\,.
\end{equation}

\medskip\noindent
\textbf{Identification of the limit.}
Recall that
\[
\mathcal Q_\kappa(z,x)
=
\Psi_\kappa(zx)+r_\kappa(z)\Phi_\kappa(zx)\,.
\]
By the regularity properties established above,
$\mathcal Q_\kappa(z,\cdot)$ belongs to
$L^2((0,1);\mathbb C^2)$ for every $z$ outside
the spectrum, including $z=0$ by holomorphic
continuation.

\medskip\noindent
For every $v\in C_c^\infty((0,1);\mathbb R^2)\,$, the end of section \ref{subsec:supcompact} combine with the identities \eqref{eq:pointwise-KS} and \eqref{eq:L2-KS-kernel-convergence} imply
\[
\int_0^1
\left[S(z,x)-\mathcal Q_\kappa(z,x)\right]
\cdot v(x)\,dx=0\,.
\]
By density,
\[
S(z,\cdot)=\mathcal Q_\kappa(z,\cdot)
\quad\text{in }L^2((0,1);\mathbb C^2)\,.
\]
Therefore, for every
$v\in L^2((0,1);\mathbb R^2)\,$,
Cauchy--Schwarz gives
\[
\int_0^1S_N(z,x)\cdot v(x)\,dx
\longrightarrow
\int_0^1\mathcal Q_\kappa(z,x)\cdot v(x)\,dx\,.
\]
Using \eqref{eq:differential-Phi} in the finite
sum on the left proves the following result.

\begin{prop}[Kneser--Sommerfeld formula for the MIT bag model]
	\label{prop:kneser-sommerfeld}
	Let $\kappa\in\mathbb N$ and
	$v\in L^2((0,1);\mathbb R^2)\,$.
	For every
	$z\in\mathbb C\setminus\{\lambda_n:n\in\mathbb Z\}\,$,
	\[
	\int_0^1
	\left[\Psi_\kappa(zx)+r_\kappa(z)\Phi_\kappa(zx)\right]
	\cdot v(x)\,dx
	=
	\lim_{N\to\infty}\sum_{n=-N}^{N}
	\frac{D_{(p,q)}\lambda_{\kappa,n}(0,0)\cdot v}
	{\lambda_n-z}\,.
	\]
\end{prop}


\section{Kernel of the differential of the spectral map}

In this section we prove that the kernel of the differential of the spectral map $\mathcal S_{0,1}$ is reduced to $\{0\}\,$:
\begin{prop}\label{prop:ker=0} Let $v=(v_1,v_2)\in L^2(0,1)\times L^2(0,1)$ satisfying
$$D_{(p,q)}\mathcal S_{0,1}\cdot v=0$$
then $v=0\,$.
\end{prop} 
\noindent
In a first subsection we use Proposition \ref{prop:kneser-sommerfeld} to transform the condition 
\begin{equation}\label{dlambda=01}D_{(p,q)}\lambda_{\kappa,n}(0,0)\cdot (v_1,v_2)=0\quad \text{for all } n\in\Z\end{equation}
into a coupled differential system linking $w_1=T_\kappa^{(1)}(v_1)$ and $w_2=T_\kappa^{(2)}(v_2)$ where $T_\kappa$ is the transformation operator defined in section \ref{sec:trans} and first introduced in \cite{Serier2006}. Then we make explicit this system in terms of $v_1$ and $v_2$ in the case $\kappa=0$ and $\kappa=1$ and prove that the two conditions $D_{(p,q)}\lambda_{0,n}(0,0)\cdot (v_1,v_2)=0$ and $D_{(p,q)}\lambda_{1,n}(0,0)\cdot (v_1,v_2)=0$ lead to $v_1=v_2=0\,$.

\subsection{Application of the Kneser--Sommerfeld formula}
In this section we assume that $v=(v_1,v_2)\in L^2(0,1)\times L^2(0,1)$ satisfies
\begin{equation}\label{dlambdav=0}D_{(p,q)}\lambda_{\kappa,n}(0,0)\cdot (v_1,v_2)=0 \quad \text{for all } n\in\Z\end{equation}
for some integer value of $\kappa$.
Then applying Proposition \ref{prop:kneser-sommerfeld} we get
$$\int_0^1
	\left[
	\Psi_\kappa(zx)+r_\kappa(z)\Phi_\kappa(zx)
	\right]\cdot v(x)\,dx
	=0 \quad\text{for all } z\in\C\,.$$
Now for $\kappa$ integer we have (recall that $\nu=\kappa+\frac12$) $ Y_\nu=(-1)^{\kappa+1}J_{-\nu}$ and thus
$$ r_\nu(z)=(-1)^\kappa \frac{J_{-\nu+1}(z)+J_{-\nu}(z)}{J_{\nu-1}(z)-J_\nu(z)}\,.$$
Therefore we get
\begin{equation}\label{dlambda=02}\int_0^1
	\left[
	(J_{\nu-1}(z)-J_\nu(z))\Psi_\kappa(zx)+(-1)^\kappa(J_{-\nu+1}(z)+J_{-\nu}(z))\Phi_\kappa(zx)
	\right]\cdot v(x)\,dx
	=0\quad \text{for all } z\in\C\,.\end{equation}

\subsection{Transformation operators}\label{sec:trans}

We recall the definition of the transformation operators introduced in the work of Serier \cite{Serier2006}.
Such operators first appeared in the seminal paper of Guillot and Ralston
\cite{GR} in connection with the inverse spectral problem for the radial 
Schr\"odinger operator (the case $\kappa=1$). They were later extended to general
integer $\kappa$ by Rundell and Sacks \cite{RuSa01}, and subsequently refined in 
\cite{Ser07}. 

\medskip
\noindent
In the AKNS setting, Serier constructed similar operators adapted to the
first-order matrix structure. A key difference with the Schr\"odinger case is
that the inverse operators have a more favorable structure.

\medskip\noindent
The next  lemma is taken from \cite{Serier2006}.
\begin{lemma}\label{akns_lemme_reduction_indice}
For each $\kappa\in\mathbb{N}\,$, define the operator
\[
S_{\kappa+1}:\;L^2(0,1)^2\longrightarrow L^2(0,1)^2\,,
\qquad
S_{\kappa+1}[p,q]:=\bigl(S_{\kappa,1}[p],\,S_{\kappa,2}[q]\bigr)\,,
\]
where
\[
S_{\kappa,1}[p](x)=p(x)-2(2\kappa+1)x^{2\kappa}\!\int_x^1\frac{p(t)}{t^{2\kappa+1}}dt\,,
\qquad
S_{\kappa,2}[q](x)=q(x)-2(2\kappa+1)x^{2\kappa+1}\!\int_x^1\frac{q(t)}{t^{2\kappa+2}}dt\,.
\]
We also set $S_0:=\mathrm{Id}\,$.
 The functions $\Phi_\kappa$ and $\Psi_\kappa$ (defined in \eqref{eq:Phikappa}) satisfy the reduction relations
\[
\Phi_{\kappa+1}=-S_{\kappa+1}^\ast[\Phi_\kappa]\,,
\qquad
\Psi_{\kappa+1}=-S_{\kappa+1}^\ast[\Psi_\kappa]\,.
\]
\end{lemma}

\noindent
We now consider the composite operator \(T_\kappa\), obtained by composing the index-reduction operators \(S_1,\dots,S_\kappa\,\), which carries Bessel kernels to trigonometric ones.

\begin{lemma}\label{akns-akns_lemme_bessel_sinus}
For every $\kappa\in\mathbb{N}\,$, define
\[
T_\kappa=(-1)^{\kappa+1}S_\kappa S_{\kappa-1}\cdots S_1\,,
\qquad
T_0:=-S_0\,.
\]
Write $T_\kappa[p,q]=(T_\kappa^1[p],\,T_\kappa^2[q])\,$. Then:
\begin{enumerate}[(i)]
\item $T_\kappa$ is bounded and injective, and for all $(p,q) \in L^2(0,1)^2$ and all $\lambda\in\mathbb{C}\,$,
\[
\int_0^1\Phi_\kappa(\lambda t)\cdot\binom{p(t)}{q(t)}\,dt
=
\int_0^1
\binom{\sin(2\lambda t)}{\cos(2\lambda t)}\cdot T_\kappa[p,q](t)\,dt\,,
\]
\[
\int_0^1\Psi_\kappa(\lambda t)\cdot\binom{p(t)}{q(t)}\,dt
=
\int_0^1
\binom{\cos(2\lambda t)}{-\sin(2\lambda t)}\cdot T_\kappa[p,q](t)\,dt\,.
\]

\item The adjoint $T_\kappa^\ast$ satisfies 
\[
\Phi_\kappa(\lambda x)
= T_\kappa^\ast
\binom{\sin(2\lambda \cdot)}{\cos(2\lambda \cdot)}(x)\,,
\qquad
\Psi_\kappa(\lambda x)
= T_\kappa^\ast
\binom{\cos(2\lambda \cdot)}{-\sin(2\lambda \cdot)}(x)\,.
\]

\end{enumerate}
\end{lemma}

\noindent
Using Lemma \ref{akns-akns_lemme_bessel_sinus} assertion (ii) we get that \eqref{dlambda=02} reads
\begin{equation}\label{dlambda=03}\int_0^1
	\left[
	(J_{\nu-1}(z)-J_\nu(z))\binom{\cos2 zx}{-\sin2 zx}+(-1)^\kappa(J_{-\nu+1}(z)+J_{-\nu}(z))\binom{\sin2 zx}{\cos2 zx}
	\right]\cdot T_\kappa [v](x)\,dx
	=0\, , \ \text{for all } z\in\C\,.\end{equation}


\vspace{0.2cm}\noindent
For later use, we recall the explicit formulas for Bessel functions of
half-integer order together with the associated polynomials introduced
in \cite[10.1.19--20]{AS64}.  
For $\kappa=0,1,2,\dots$ and $z\in\C$, one has the classical representations
\begin{align}
J_{\kappa+\tfrac12}(z)
&=
\sqrt{\frac{2}{\pi z}}\,
\Bigl(
   P_\kappa\!\left(\tfrac1z\right)\,\sin z
   \;-\;
   Q_{\kappa-1}\!\left(\tfrac1z\right)\,\cos z
\Bigr)\,,
\label{besseldemientier1}
\\[0.3em]
J_{-\kappa-\tfrac12}(z)
&=
(-1)^{\kappa}\sqrt{\frac{2}{\pi z}}\,
\Bigl(
   P_\kappa\!\left(\tfrac1z\right)\,\cos z
   \;+\;
   Q_{\kappa-1}\!\left(\tfrac1z\right)\,\sin z
\Bigr)\,.
\label{besseldemientier2}
\end{align}

\vspace{0.15cm}\noindent
The polynomials \(P_\kappa\) and \(Q_\kappa\), each of degree~\(\kappa\), are generated by the three-term recurrences 
\begin{align*}
P_{\kappa+1}(t)
&=(2\kappa+1)\,t\,P_\kappa(t)-P_{\kappa-1}(t)\,,
\qquad \kappa\ge1\,,
\\
Q_{\kappa+1}(t)
&=(2\kappa+3)\,t\,Q_\kappa(t)-Q_{\kappa-1}(t)\,,
\qquad \kappa\ge0\,,
\end{align*}
with initial values
\begin{equation}\label{P0}
P_0(t)=1\,,\qquad P_1(t)=t\,,
\qquad
Q_{-1}(t)=0\,,\qquad Q_0(t)=1\,.
\end{equation}
\noindent
Finally we introduce the polynomials
\begin{align*}
U_{\kappa}(t)
&=P_\kappa(t)-P_{\kappa-1}(t)\,,
\qquad \kappa\ge1\,,
\\
V_{\kappa}(t)
&=Q_{\kappa-1}(t)-Q_{\kappa-2}(t)\,,
\qquad \kappa\ge1\,.
\end{align*}
Using \eqref{besseldemientier1} and \eqref{besseldemientier2} and the definition of $U_\kappa$ and $ V_\kappa$, \eqref{dlambda=03} can be rewritten as
\begin{equation}\label{dlambda=04}\int_0^1
	\left[
	V_\kappa(\frac1z)\binom{\cos z(2x-1)}{-\sin z(2x-1)}+U_\kappa(\frac1z) \binom{\sin z(2x-1)}{\cos z(2x-1)}
	\right]\cdot T_\kappa v(x)\,dx
	=0\quad \text{for all } z\in\C\,.\end{equation}
We now introduce the polynomials $\tilde U_\kappa(z)=z^\kappa U_\kappa(\frac1z)$ and $\tilde V_\kappa(z)=z^\kappa V_\kappa(\frac1z)$ that we decompose in even and odd part: $\tilde U_\kappa=\tilde U^{\text{even}}_\kappa+\tilde U^{\text{odd}}_\kappa$, $\tilde V_\kappa=\tilde V^{\text{even}}_\kappa+\tilde V^{\text{odd}}_\kappa$. We notice that
\begin{align*}
z^\kappa U_\kappa(\frac1z) \binom{\sin z(2x-1)}{\cos z(2x-1)}&=  \tilde U^{\text{even}}_\kappa(\tilde D)\binom{\sin z(2x-1)}{\cos z(2x-1)}-i \tilde U^{\text{odd}}_\kappa(\tilde D) \binom{\cos z(2x-1)}{-\sin z(2x-1)}\\
z^\kappa V_\kappa(\frac1z) \binom{\cos z(2x-1)}{-\sin z(2x-1)}&=  \tilde V^{\text{even}}_\kappa(\tilde D)\binom{\cos z(2x-1)}{-\sin z(2x-1)}+i \tilde V^{\text{odd}}_\kappa(\tilde D) \binom{\sin z(2x-1)}{\cos z(2x-1)}\ ,
\end{align*}
where $\tilde D=\frac1{2i}\frac d{dx}\,$. Therefore introducing (in the distributional sense)
\begin{align*}
f_\nu  &=\big( \tilde U^{\text{even}}_\kappa(\tilde D)-i \tilde V^{\text{odd}}_\kappa(\tilde D)\big) T_\kappa(v)\,,\\
g_\nu &=\big( \tilde V^{\text{even}}_\kappa(\tilde D)+i \tilde U^{\text{odd}}_\kappa(\tilde D)\big) T_\kappa(v)\,,
\end{align*}
\eqref{dlambda=04} can be rewritten as
$$\int_0^1 \big(\sin (2x-1)z\ (f_{\nu1}(x)-g_{\nu2}(x))+\cos (2x-1)z\ (f_{\nu2}(x)+g_{\nu1}(x))\big)dx=0 \quad \text{for all }z\in\C\setminus\{0\}\,.$$
Then changing $z$ in $-z$ and taking the sum and the difference we get for all $z\in\C\setminus\{0\}$
$$\int_0^1\cos z(2x-1)\ (f_{\nu2}+g_{\nu1})(x) dx=0 \text{ and }\int_0^1\sin z(2x-1)\ (f_{\nu1}-g_{\nu2})(x) dx=0\,,$$
which finally leads to
\begin{equation}\label{fg}
\begin{cases}
f_{\nu1}-g_{\nu2}\quad  \text{even with respect to }1/2\,,\\
f_{\nu2}+g_{\nu1}\quad  \text{odd with respect to }1/2\,.
\end{cases}
\end{equation}

\subsection{The cases $\kappa=0$ and $\kappa=1$}

In this section we reinterpret the assumption \eqref{dlambda=01} in terms of $v$ in the two cases $\kappa=0$ and $\kappa=1$.\\
When $\kappa=0$, direct application of \eqref{dlambda=02} using \eqref{eq:HalfInt}, \eqref{eq:Phikappa} and $Y_\nu=(-1)^{\kappa+1}J_{-\nu}$ leads to 
$$\int_0^1 \big(\sin (z(2x-1)+\frac\pi4)v_1(x)+\cos (z(2x-1)+\frac\pi4)v_2(x)\big)dx=0 \quad \text{for all }z\in\C\,.$$
Then changing $z$ in $-z$ and taking the sum and the difference we get
$$\int_0^1\cos z(2x-1)\ (v_1+v_2)(x) dx=0 \text{ and }\int_0^1\sin z(2x-1)\ (v_1-v_2)(x) dx=0\,,$$
which leads to $v_1+v_2$ odd and $v_1-v_2$ even and thus
\begin{equation}\label{kappa=0} v_2(x)=-v_1(1-x).\end{equation}
In the case $\kappa=1$, using \eqref{P0}, we first compute $\tilde U_1(t)=1-t$ and $\tilde V_1(t)=t$ and thus
$f_1=(1-i\tilde D)T_1(v)$ and $g_1=-i\tilde D T_1(v)\,$. Thus, denoting $w=T_1(v)$, \eqref{fg} yield
\begin{equation}\label{w}
\begin{cases}
-w_1'+w_2'+2w_1 \quad \text{even},\\
w_1'+w_2'-2w_2 \quad \text{odd}.
\end{cases}
\end{equation}
which in turn implies
\begin{equation}\label{w'}
\begin{cases}
-w_1''+w_2''+2w'_1 \quad \text{odd},\\
w_1'''+w_2'''-2w''_2 \quad \text{odd}.
\end{cases}
\end{equation}
On the other hand,  by Lemma \ref{akns_lemme_reduction_indice} 
\begin{equation}\label{wv}
\begin{cases}
w_1=v_1-2\int_x^1 \frac{v_1(t)}{t}\,dt\,,\\
w_2=v_2-2x\int_x^1 \frac{v_2(t)}{t^2}\,dt\,.
\end{cases}
\end{equation}
Inserting \eqref{wv} in \eqref{w'}we get
\begin{equation}\label{v}
\begin{cases}
-v_1''+v_2''+\frac2x(v_2'-v_1')+2v_1'+\big(\frac2{x^2}+\frac4x\big)v_1  \quad \text{odd},\\
v_1'''+v_2'''+\frac2x(v_1''+v_2'')-2v_2''-\frac4{x^2}v_1'+\big(-\frac2{x^2}-\frac4x\big)v_2'+\frac4{x^3}v_1 \quad \text{odd}.
\end{cases}
\end{equation}

\subsection{Kernel of the differential of the spectral map $\mathcal S_{0,1}$}
In this section we assume that $v=(v_1,v_2)\in L^2(0,1)\times L^2(0,1)$ is in the kernel of $D_{(p,q)}\mathcal S_{0,1}(0,0)\,$, i.e. that
 \eqref{dlambda=01} for both $\kappa=0$ and $\kappa=1\,$. We aim at proving $v=0$.\\
First we decompose $v_1$  in an even part, denoted  $f$, and an odd part denoted $g$ : 
\begin{equation}\label{vfg}v_1=f+g\,,\mbox{ with } f\text{ even and }g \text{ odd}.\end{equation}
 By \eqref{kappa=0} we get $v_2=-f+g$ and then \eqref{v} leads to
the following system
\begin{equation}\label{benoitsyst}
\begin{cases}
	2g'''(x)+2f''(x)+2B_1(x)g''(x)
	-\bigl(3A_2(x)+2A_1(x)\bigr)g'(x)\\
	\qquad
	+\bigl(2B_1(x)-B_2(x)\bigr)f'(x)
	+2A_3(x)f(x)+2B_3(x)g(x)=0\,,\\[0.6em]
	-2f''(x)-2B_1(x)f'(x)+2g'(x)
	+\bigl(2A_1(x)+A_2(x)\bigr)f(x)
	+\bigl(2B_1(x)+B_2(x)\bigr)g(x)=0\,,
\end{cases}
\end{equation}
 where
\[
A_j(x)=\frac{1}{x^j}+\frac{1}{(1-x)^j}\, \mbox{ and } 
B_j(x)=\frac{1}{x^j}-\frac{1}{(1-x)^j},\quad j=1,2,3\,.
\]
\subsubsection{Indicial equation at \(x=0\)}

As \(x\to 0\)\,, we have \(A_j(x)\sim B_j(x)\sim x^{-j}\)\,. Assume that
\[
f(x)\sim a x^\alpha,\qquad g(x)\sim b x^\beta\,,
\qquad a,b\ne 0\,.
\]
The leading terms of the two equations $E_1$, $E_2$ in \eqref{benoitsyst} are
\[
\begin{aligned}
	E_1 &\sim
	a(2-\alpha)x^{\alpha-3}
	+b(2\beta^3-4\beta^2-\beta+2)x^{\beta-3}\,,\\
	E_2 &\sim
	a(1-2\alpha^2)x^{\alpha-2}
	+b x^{\beta-2}\,.
\end{aligned}
\]
If \(\alpha<\beta\)\,, cancellation would require both
\(\alpha=2\) and \(\alpha=\pm 1/\sqrt{2}\)\,, which is impossible.
If \(\beta<\alpha\)\,, the leading term \(b x^{\beta-2}\) in \(E_2\)
cannot vanish. Hence a nontrivial power-law mode must have
\(\alpha=\beta=\rho\)\,.

\medskip\noindent
The indicial equation is therefore
\[
\det
\begin{pmatrix}
	2-\rho & 2\rho^3-4\rho^2-\rho+2\\
	1-2\rho^2 & 1
\end{pmatrix}=0\,,
\]
or equivalently
\[
4\rho^2(\rho+1)(\rho-1)(\rho-2)=0\,.
\]
Thus the indicial exponents are \(-1,0,0,1,2\) and the Frobenius recurrence must be
checked for a possible logarithmic solution. Any term proportional
to \(\log x\) is still square-integrable near \(0\). The exponent
\(\rho=-1\)\,, however, gives a non-square-integrable mode.

\subsubsection{Equivalent first-order system}

Set
\begin{equation}\label{hkr}
h:=-f',\qquad k:=-g'\,,\qquad r:=g''\,,
\end{equation}
and define
\[
\alpha:=A_1+\frac{A_2}{2},\qquad
\beta:=B_1+\frac{B_2}{2},\qquad
\gamma:=\frac{3A_2+2A_1-2}{2}\,.
\]
Then the system \eqref{benoitsyst} is equivalent on \((0,1)\) to
\begin{equation}\label{eq:system}
	\begin{cases}
		f'=-h\,,\\
		g'=-k\,,\\
		h'=-\alpha f-\beta g-B_1h+k\,,\\
		k'=-r\,,\\
		r'=-(\alpha+A_3)f-(\beta+B_3)g-\dfrac{B_2}{2}h-\gamma k-B_1r\,.
	\end{cases}
\end{equation}

\subsubsection{End of the proof of Proposition \ref{prop:ker=0}}
Assume that $D_{(p,q)}\mathcal S_{0,1}\cdot v=0$. Then by construction $f,g,h,k,r$ defined by \eqref{vfg} and \eqref{hkr} satisfy the differential system \eqref{eq:system}. By parity we know that \(f'(\tfrac12)=g(\tfrac12)=g''(\tfrac12)=0\) but furthermore, going back to \eqref{w} we know that $w_1'(\tfrac12)+w_2'(\tfrac12)-2w_2(\tfrac12)=0$. After calculation this last relation leads to $g'(\frac12)+f(\frac12)=0\,$. 

\noindent
As a consequence if \(f(\tfrac12)=0\), then all the
Cauchy data in \eqref{eq:system} at $\tfrac12$ are zero, and uniqueness gives
\(f=g=0\). 

\noindent
It remains to consider the case \(f(\tfrac12)\ne0\) and by linearity we may
normalize \(f(\tfrac12)=1\). By parity, the Cauchy data at \(x=\tfrac12\) are
therefore 
\begin{equation}\label{eq:initial-data}
	(f,g,h,k,r)\big|_{x=1/2}=(1,0,0,1,0)\,.
\end{equation}

\medskip\noindent
Define
\begin{equation}\label{eq:Q}
	Q:=\alpha f+\beta g+B_1h-k\,.
\end{equation}
The third equation of \eqref{eq:system} gives
\begin{equation}\label{eq:hQ}
	h'=-Q\,.
\end{equation}
We also set
\begin{equation}\label{eq:D}
	D:=-Q'\,.
\end{equation}
The coefficient functions satisfy
\[
A_j'=-jB_{j+1}\,,\qquad B_j'=-jA_{j+1}\,.
\]
We  differentiate \eqref{eq:Q} and substitute
\eqref{eq:system}, we obtain
\begin{align}
	Q'
	&=r-(B_2+B_3+B_1\alpha)f
	-(A_2+A_3+B_1\beta)g \notag\\
	&\qquad-(\alpha+A_2+B_1^2)h-\frac{B_2}{2}k\,.
	\label{eq:Qprime}
\end{align}
Using \(D=-Q'\)\,, and differentiating once more with Mathematica, we get
\begin{equation}\label{eq:Dprime}
	D'=-\left(
	c_Ff+c_Gg+c_Hh+c_Kk+\frac{B_2}{2}r+c_QQ+B_1D
	\right)\,,
\end{equation}
where
\begin{align}
	c_F&:=A_3+3A_4-\alpha\,,\label{eq:cF}\\
	c_G&:=B_3+3B_4-\beta\,,\label{eq:cG}\\
	c_H&:=\frac32B_2+4B_3\,,\label{eq:cH}\\
	c_K&:=2A_2+2A_3-\gamma\,,\label{eq:cK}\\
	c_Q&:=\alpha+2A_2\,.\label{eq:cQ}
\end{align}
Thus, using \eqref{eq:system}, \eqref{eq:hQ}, and \eqref{eq:Dprime}, we obtain 
 the following system:
\begin{equation}\label{eq:augmented-system}
	\begin{cases}
		f'=-h\,,\\
		g'=-k\,,\\
		h'=-Q\,,\\
		k'=-r\,,\\
		r'=- \Big(  (\alpha+A_3)f+(\beta+B_3)g
		+\frac{B_2}{2}h+\gamma k+B_1r            \big)\,,\\
		Q'=-D\,,\\
		D'=-\big(c_Ff+c_Gg+c_Hh+c_Kk+\dfrac{B_2}{2}r+c_QQ+B_1D \big)\,.
	\end{cases}
	\end{equation}
At \(x:=\tfrac12\)\,,  we get the Cauchy data:
\[
(f,g,h,k,r,Q,D)(\tfrac12)=(1,0,0,1,0,7,0)\,.
\]

\paragraph{Signs of the coefficients.}
On \(0<x\leq\tfrac12\,\), we have
\[
A_1\leq A_2\leq A_3\leq A_4\,,\qquad
0\leq B_1\leq B_2\leq B_3\leq B_4\,.
\]
It follows from \eqref{eq:cF}--\eqref{eq:cQ} that all five coefficients
\(c_F,c_G,c_H,c_K,c_Q\) are nonnegative. Indeed,
\[
c_F\geq \frac52 A_3>0\,,\qquad
c_G\geq \frac52 B_3\geq0\,,\qquad
c_H\geq0\,,\qquad c_Q>0\,,
\]
and
\[
c_K=\frac12A_2+2A_3-A_1+1>0\,.
\]
The coefficients in the fifth equation in \eqref{eq:augmented-system} are also nonpositive on \((0,\tfrac12]\)\,, with
\(\alpha+A_3>0\) and \(\gamma>0\)\,.

\medskip\noindent
Moreover, \eqref{eq:augmented-system} gives
\[
g'(\tfrac12)=-1\,,\qquad h'(\tfrac12)=-7\,,\qquad r'(\tfrac12)=-39\,.
\]
At \(\tfrac12\)\,, the coefficients in \eqref{eq:cF}--\eqref{eq:cQ} satisfy
\(c_F=104\)\,, \(c_K=33\)\,, and \(c_Q=24\)\,. Thus
\[
D'(\tfrac12)=-(104+33+24\cdot7)=-305\,.
\]
It follows that all seven components are nonnegative on an interval in the form 
\([\tfrac12-\delta,\tfrac12]\) for small enough  $\delta>0$\,.

\medskip\noindent
Let
\[
x_0:=\inf\left\{x\in(0,\frac12]:
f,g,h,k,r,Q,D\geq0\ \text{on }[x,\frac12]\right\}\,.
\]
By the previous discussion, the set is nonempty and we have
\[
x_0\leq \frac12-\delta\,.
\]
We now argue by contradiction and assume that \(x_0>0\). By continuity, all seven functions are nonnegative on
\([x_0,\frac12]\). On this interval, \eqref{eq:augmented-system} shows that
each is nonincreasing as a function of \(x\). Therefore, for $x \in [x_0,\frac12]\,$,
\begin{align}
	Q(x)&\geq Q(\frac{1}{2}) = 7\,,\label{eq:Q-lower}\\
	h(x)&=\int_x^{\frac12}Q(s)\,ds
	\geq7\,(\frac12-x)\,,\label{eq:h-lower}\\
	f(x)&=1+\int_x^{\frac12}h(s)\,ds
	\geq1+\frac72\,(\frac12-x)^2\,,\label{eq:f-lower}\\
	k(x)&\geq k(\frac{1}{2}) = 1\,,\label{eq:k-lower}\\
	g(x)&=\int_x^{\frac12}k(s)\,ds\geq \frac12-x\,.\label{eq:g-lower}
\end{align}
We also have  on the interval $[x_0, \frac12]$\,,
$$
(\alpha+A_3)f+(\beta+B_3)g
+\frac{B_2}{2}h+\gamma k+B_1r\geq(\alpha+A_3)f>0
$$
 and 
 $$
 c_Ff+c_Gg+c_Hh+c_Kk+\dfrac{B_2}{2}r+c_QQ+B_1D\geq c_Ff+c_QQ>0 \,.
 $$
Using \eqref{eq:augmented-system} and integrating $r'$ and $D'$ from \(x_0\) to \(\frac12\) we show that
\(r(x_0)>0\) and \(D(x_0)>0\)\,. In particular, all seven
components are strictly positive at \(x_0\). By continuity, they remain
positive on some interval to the left of \(x_0\), contradicting the
definition of \(x_0\). 

\medskip\noindent
Thus \(x_0=0\), and all seven components are
nonnegative on \((0,\frac12]\)\,. In particular, 
\eqref{eq:f-lower} and \eqref{eq:g-lower} imply
\begin{equation}\label{eq:strict-f}
f(x)>1\,, \qquad g(x)>0\,, 
	\qquad 0<x<\frac12\,.
\end{equation}

\medskip\noindent
It remains to use the Frobenius expansion at $x=0\,$. The indicial roots
are $-1,0,0,1,2$\,. Thus, as $x \to 0^+$\,, 
\[
f(x)=\frac{A}{x}+B\log x+C+o(1)\,,\qquad
g(x)=\frac{A'}{x}+B'\log x+C'+o(1)\,.
\]
Since $f,g\in L^2(0,1/2)$, we have $A=A'=0$. Substituting the
Frobenius expansions into the second equation of \eqref{benoitsyst}
and multiplying by $x^2$, we obtain
\[
(B+B')\log x+C+C'=o(1),
\]
since the constant terms from $x^2f''$ and $x^2B_1f'$ cancel.
Thus $B'=-B$ and $C'=-C$, so $f(x)+g(x)\to0$ as $x\to0^+$.
This contradicts \eqref{eq:strict-f}, which gives $f(x)+g(x)>1$.
Therefore $f=g=0$, and hence $v_1=v_2=0$. This completes the proof
of Proposition~\ref{prop:ker=0}.


\section{Invertibility of the linearized spectral map}
\label{sec:mit-closed-range}

In this section, we prove that the spectral differential at the zero
potential is a bounded isomorphism for $(\kappa_1,\kappa_2)=(0,1)\,$. More precisely, we show that it is Fredholm of index zero and we use the injectivity
established in the preceding section.
As a consequence, the analytic inverse function theorem implies that
the spectral map is a local real-analytic diffeomorphism near zero.

\medskip\noindent
We follow the approach of \cite[Section~7]{GGHN2026}.
Although the MIT bag spectrum is not symmetric about zero,
a hidden symmetry allows us to split the space of pairs of functions
into two orthogonal subspaces.
The proof strategy is explained in Subsection~6.2.

\subsection{Preliminaries}

Set
\[
\mathcal H=L^2((0,1);\mathbb R)\, ,\qquad
X=\mathcal H\times\mathcal H \, ,\qquad
Y=\ell^2(\mathbb Z;\mathbb R)\times\ell^2(\mathbb Z;\mathbb R) \, .
\]
The spectral differential at the zero potential is the operator
\[
S:=D_{(p,q)}\mathcal S_{0,1}(0,0):X\longrightarrow Y\, ,
\]
given, for $v=(v_1,v_2)\in X\,$, by
\[
Sv=
\left(
\bigl(D_{(p,q)}\lambda_{0,n}(0,0)\cdot v\bigr)_{n\in\mathbb Z},
\bigl(D_{(p,q)}\lambda_{1,n}(0,0)\cdot v\bigr)_{n\in\mathbb Z}
\right).
\]
We aim at proving
\begin{prop}\label{iso}
$S:=D_{(p,q)}\mathcal S_{0,1}(0,0)$ is a bounded isomorphism from $X$ onto $Y$.
\end{prop}

\noindent
Writing $\lambda_{\kappa,n}:=\lambda_{\kappa,n}(0,0)\,$,
we recall from \eqref{eq:differential-Phi} that
\begin{equation}\label{eq:mit-cr-differential}
	D_{(p,q)}\lambda_{\kappa,n}(0,0)\cdot v
	=a_{\kappa,n}\int_0^1
	\Phi_\kappa(\lambda_{\kappa,n}x)\cdot v(x)\,dx\, ,
	\qquad
	a_{\kappa,n}:=\frac{2c_{\kappa,n}^2}{\pi} >0 \, .
\end{equation}
We first compute the normalization factors explicitly
for $\kappa=0$ and $\kappa=1\,$.

\begin{lemma}
	For every $n\in\mathbb Z\, $, the normalization factors satisfy
	\[
	a_{0,n}=1\, ,
	\qquad
	a_{1,n}
	=1+\frac{1}{2\lambda_{1,n}(\lambda_{1,n}-1)} \,.
	\]
	In particular, 
	\[
	a_{1,n}=1+O(|n|^{-2})
	\quad , \quad |n|\to\infty\, .
	\]
\end{lemma}

\begin{proof}
	By \eqref{eq:mit-cr-differential},
	\eqref{eq:c-normalization} and \eqref{eq:characteristic},
	\begin{equation}\label{eq:mit-a-normalization}
		a_{\kappa,n}
		=\frac{1}{
			\pi(\lambda_{\kappa,n}-\kappa)
			J_{\kappa-\frac12}(\lambda_{\kappa,n})^2}\, .
	\end{equation}
    Formula \eqref{eq:HalfInt} gives
	\[
	J_{-\frac12}(\lambda)
	=\sqrt{\frac{2}{\pi\lambda}}\cos\lambda \, .
	\]
	Since $\lambda_{0,n}=\pi n+\pi/4$\,, as recalled in
	\eqref{eq:l0-eigs}, we obtain
	\begin{equation}\label{eq:mit-a-zero}
		a_{0,n}
		=\frac{1}{2\cos^2\lambda_{0,n}}
		=1 \, .
	\end{equation}
For $\kappa=1$, write $\lambda=\lambda_{1,n}$ and recall that $\lambda\ne0$.
	Formula \eqref{besseldemientier1} gives
	\[
	J_{\frac12}(\lambda)
	=\sqrt{\frac{2}{\pi\lambda}}\sin\lambda\,,
	\qquad
	J_{\frac32}(\lambda)
	=\sqrt{\frac{2}{\pi\lambda}}
	\left(\frac{\sin\lambda}{\lambda}-\cos\lambda\right)\,.
	\]
	The characteristic equation \eqref{eq:characteristic}
	therefore yields
	\begin{equation}\label{eq:mit-characteristic-kappa-one-trig}
		\cos\lambda
		=\left(\frac1\lambda-1\right)\sin\lambda\,.
	\end{equation}
	In particular, $\sin\lambda\ne0$, and hence
	\[
	\frac1{\sin^2\lambda}
	=1+\left(1-\frac1\lambda\right)^2.
	\]
	Substituting into \eqref{eq:mit-a-normalization} gives
	\begin{equation}\label{eq:mit-a-one-exact}
		a_{1,n}
		=\frac{\lambda}{2(\lambda-1)}
		\left[1+\left(1-\frac1\lambda\right)^2\right]
		=1+\frac{1}{2\lambda(\lambda-1)}\,.
	\end{equation}
	Finally, \eqref{eq:AsymptAKNS} gives
	$\lambda_{1,n}=\pi n+O(1)$ as $|n|\to\infty$.
	Together with \eqref{eq:mit-a-one-exact}, this implies
	\begin{equation}\label{eq:mit-a-one-asymptotics}
		a_{1,n}=1+O(|n|^{-2})\,,
		\mbox{ as } |n|\to\infty\,.
	\end{equation}
\end{proof}

\medskip\noindent
Thus multiplication by $(a_{\kappa,n})$ defines a bounded isomorphism
$\mathcal D:Y\to Y\,$. Consequently, $S:X\to Y$ is a bounded
isomorphism if and only if the normalized operator
\[
L:=\mathcal D^{-1}S:X\to Y
\]
is a bounded isomorphism.

\medskip\noindent
For later use, introduce the bounded\footnote{The fact that $H$ and $L$ are continuous from $\mathcal{H}$ to $\mathcal{H}$ easily follows from  Hardy's inequality. } integral operators on $\mathcal{H}$
\begin{equation}\label{HJ}
	(Hf)(x):=\int_x^1\frac{f(t)}{t}\,dt\,,
	\qquad
	(Jf)(x):=x\int_x^1\frac{f(t)}{t^2}\,dt\,,
\end{equation}
and set
\begin{equation}\label{AB}
	A:=I-2H=S_{0,1}\,,\qquad B:=I-2J=S_{0,2}\,,\qquad
	\mathcal A(p,q):=(Ap,Bq)\,.
\end{equation}
The operators $S_{0,1}$ and $S_{0,2}$ in
\eqref{AB} are the scalar transformation operators
defined in Lemma~\ref{akns_lemme_reduction_indice} with $\kappa=0$. By \eqref{eq:mit-cr-differential}, the normalized differential is
\begin{equation}\label{eq:mit-cr-L-Phi}
	Lv=
	\left(
	\left(\int_0^1
	\Phi_0(\lambda_{0,n}x)\cdot v(x)\,dx\right)_{n\in\mathbb Z},
	\left(\int_0^1
	\Phi_1(\lambda_{1,n}x)\cdot v(x)\,dx\right)_{n\in\mathbb Z}
	\right)\,.
\end{equation}
By Lemma~\ref{akns-akns_lemme_bessel_sinus},
the transformation operators satisfy
\[
T_0v {=-S_0v} = -v\,,
\qquad
T_1v=S_1v=(S_{0,1}p,S_{0,2}q)=(Ap,Bq)\,,
\]
where $v=(p,q)$. Applying Lemma~\ref{akns-akns_lemme_bessel_sinus}\,(i) to
\eqref{eq:mit-cr-L-Phi}, we obtain
\begin{align}
	(Lv)_{0,n}
	&=-\int_0^1
	\bigl[\sin(2\lambda_{0,n}x)p(x)
	+\cos(2\lambda_{0,n}x)q(x)\bigr]\,dx\,,
	\label{eq:mit-cr-kappa-zero}\\
	(Lv)_{1,n}
	&=\int_0^1
	\bigl[\sin(2\lambda_{1,n}x)(Ap)(x)
	+\cos(2\lambda_{1,n}x)(Bq)(x)\bigr]\,dx\,.
	\label{eq:mit-cr-kappa-one}
\end{align}
Here $(Lv)_{0,n}$ and $(Lv)_{1,n}$ denote the $n$th entries
of the first and second components of $Lv$\,, respectively.

\subsection{Proof strategy for Proposition \ref{iso}}

Recall that $2\lambda_{0,n}=2\pi n+\pi/2$\,, while
\eqref{eq:AsymptAKNS} gives
\[
2\lambda_{1,n}=
\begin{cases}
	2\pi n+3\pi/2+\varepsilon_n\,,& n\ge0\,,\\
	2\pi n-\pi/2+\varepsilon_n\,,& n\le-1\,,
\end{cases}
\quad \mbox{ with } 
(\varepsilon_n)_{n\in\mathbb Z}\in\ell^2(\mathbb Z)\,.
\]
Replacing the eigenvalues in
\eqref{eq:mit-cr-kappa-zero}--\eqref{eq:mit-cr-kappa-one}
by these leading terms suggests introducing the trigonometric vectors
\[
e_n^+(x)=
\begin{pmatrix}
	\sin((2\pi n+\pi/2)x)\\
	\cos((2\pi n+\pi/2)x)
\end{pmatrix}\,,
\qquad
e_n^-(x)=
\begin{pmatrix}
	\sin((2\pi n+3\pi/2)x)\\
	\cos((2\pi n+3\pi/2)x)
\end{pmatrix}\,,
\qquad n\in\mathbb Z\,.
\]
For $\kappa=0\,$, all vectors $e_n^+$ occur.
For $\kappa=1\,$, the leading terms give $e_n^-$ for $n\ge0$
and $e_{n-1}^-$ for $n\le-1\,$, since
\[
2\pi n-\pi/2=2\pi(n-1)+3\pi/2\,.
\]
Thus only the mode $e_{-1}^-$ is missing. We first keep all modes and define the model operator
\[
\widehat L:X\longrightarrow Y,\qquad
\widehat Lv=
\left(
\bigl(-\langle v,e_n^+\rangle_X\bigr)_{n\in\mathbb Z},
\bigl(\langle\mathcal Av,e_n^-\rangle_X\bigr)_{n\in\mathbb Z}
\right).
\]
Although the MIT bag spectrum is not symmetric {about zero},
these trigonometric vectors have a hidden symmetry under reflection.
To see this, set $Rf(x)=f(1-x)$, so that $R^2=I$, and observe that
\[
\begin{aligned}
	\sin((2\pi n+\pi/2)(1-x))
	&=\cos((2\pi n+\pi/2)x)\,,\\
	\sin((2\pi n+3\pi/2)(1-x))
	&=-\cos((2\pi n+3\pi/2)x)\,.
\end{aligned}
\]
Hence $e_n^+\in E_+$ and $e_n^-\in E_-$, where
\[
E_+=\{(f,Rf):f\in\mathcal H\}\,,
\qquad
E_-=\{(f,-Rf):f\in\mathcal H\}\,.
\]
We show that the families $(e_n^+)_{n\in\mathbb Z}$ and
$(e_n^-)_{n\in\mathbb Z}$ are orthonormal bases of
$E_+$ and $E_-\,$, respectively, and that
$X=E_+\oplus E_-$ is an orthogonal decomposition.
This decomposition allows us to write $\widehat L$ in triangular
block form and reduce the problem to a scalar operator.
We then show that this scalar operator, and hence $\widehat L\,$,
is Fredholm of index $-1\,$.

\medskip\noindent
Finally, deleting the coordinate corresponding to $e_{-1}^-$
increases the index by one.
After reindexing, the resulting operator differs from $L$
by a compact operator, as in \cite[Section~7]{GGHN2026}.
Thus $L$ is Fredholm of index zero.
Together with the injectivity proved in the preceding section,
this shows that $L$, and hence $S$, is a bounded isomorphism.

\subsection{Trigonometric model}

We recall that  $R:\mathcal H\to\mathcal H$ denotes the reflection
\[
(Rf)(x):=f(1-x)\,.
\]
and
\[
E_+:=\{(f,Rf):f\in\mathcal H\}\,,\qquad
E_-:=\{(f,-Rf):f\in\mathcal H\}\,.
\]
These subspaces are closed and orthogonal, and $X=E_+\oplus E_-\,$.
The associated orthogonal projections are
\begin{equation}
	P_+(p,q)=\frac12(p+Rq,q+Rp)\,,\qquad
	P_-(p,q)=\frac12(p-Rq,q-Rp)\,.
	\label{eq:mit-cr-projections}
\end{equation}

\begin{lemma}
	\label{lem:mit-cr-bases}
	The families
	\[
	e_n^+(x):=
	\begin{pmatrix}
		\sin((2\pi n+\pi/2)x)\\ \cos((2\pi n+\pi/2)x)
	\end{pmatrix},\qquad
	e_n^-(x):=
	\begin{pmatrix}
		\sin((2\pi n+3\pi/2)x)\\ \cos((2\pi n+3\pi/2)x)
	\end{pmatrix}\,,\quad \mbox{ for } n\in\mathbb Z\,,
	\]
	are orthonormal bases of $E_+$ and $E_-$, respectively.
\end{lemma}

\begin{proof}
We have
	\[
	\langle e_n^\pm,e_m^\pm\rangle_X
	=\int_0^1\cos(2\pi(n-m)x)\,dx
	=\delta_{nm} \, .
	\]
Hence each family is orthonormal. It remains to show that they are
complete in $E_+$ and $E_- \,$, respectively.
	The maps
	\[
	U_\pm:\mathcal H\longrightarrow E_\pm \, ,
	\qquad
	U_\pm f=\frac{1}{\sqrt{2}}(f,\pm Rf)\, ,
	\]
	are onto by the definition of $E_\pm \,$. Since $R$ preserves
	the norm of $\mathcal H$, we have
	\[
	\|U_\pm f\|_X^2
	=\frac12\bigl(\|f\|_{\mathcal H}^2+\|Rf\|_{\mathcal H}^2\bigr)
	=\|f\|_{\mathcal H}^2 \,.
	\]
	Thus $U_\pm$ are unitary. It therefore suffices to prove that
	the inverse images of our vectors form complete families
	in $\mathcal H$. Let
	\[
	s_j(x)=\sqrt{2}\sin((j+\tfrac12)\pi x) \,,
	\qquad j=0,1,2,\ldots.
	\]
The functions $s_j$ are the normalized eigenfunctions of the
Sturm--Liouville problem
\[
-u''=\lambda u\, ,\qquad u(0)=0\, ,\qquad u'(1)=0 \,,
\]
and hence form an orthonormal basis of $\mathcal H$.
For every $k\geq0$, we have
\[
U_+s_{2k}=e_k^+,\qquad
U_+s_{2k+1}=-e_{-k-1}^+,
\]
and
\[
U_-s_{2k+1}=e_k^-,\qquad
U_-s_{2k}=-e_{-k-1}^-.
\]
Thus each family $(e_n^\pm)_{n\in\mathbb Z}$ is, up to signs
and reordering, the image of the basis $(s_j)_{j\geq0}$ under
the unitary map $U_\pm$. This proves completeness in $E_\pm$.
\end{proof}

\medskip\noindent
We recall that  the model operator for $L$ is given by :
\begin{equation}
	\widehat L v
	:=\left(
	(-\langle v,e_n^+\rangle_X)_{n\in\mathbb Z},
	(\langle\mathcal A v,e_n^-\rangle_X)_{n\in\mathbb Z}
	\right).
	\label{eq:mit-cr-full-model}
\end{equation}
By Lemma~\ref{lem:mit-cr-bases}, the two blocks give the
coefficients of $-P_+v$ and $P_-\mathcal A v$ in orthonormal
bases of $E_+$ and $E_-$, respectively. In particular,
Parseval's identity gives
\begin{equation}
	\|\widehat L v\|_Y^2
	=\|P_+v\|_X^2+\|P_-\mathcal A v\|_X^2.
	\label{eq:mit-cr-parseval}
\end{equation}
Writing $v=v_++v_-$ according to $X=E_+\oplus E_-$, we have
\[
P_+v=v_+\,,
\qquad
P_-\mathcal A v=Cv_++Dv_-\,,
\]
where
\[
C=P_-\mathcal A|_{E_+}:E_+\longrightarrow E_-\,,
\qquad
D=P_-\mathcal A|_{E_-}:E_-\longrightarrow E_-\,.
\]
Thus, using these orthonormal bases, we identify $\widehat L$ with
\[
\begin{pmatrix}
	-I_{E_+} & 0\\
	C & D
\end{pmatrix}
=
\begin{pmatrix}
	-I_{E_+} & 0\\
	C & I_{E_-}
\end{pmatrix}
\begin{pmatrix}
	I_{E_+} & 0\\
	0 & D
\end{pmatrix}.
\]
The first factor is bounded and is its own inverse.
Therefore, $\widehat L$ is Fredholm if and only if $D$ is Fredholm,
and in that case
\[
\operatorname{ind}\widehat L=\operatorname{ind}D.
\]

\medskip\noindent
We now express $D$ as a scalar operator on $\mathcal H$.
\begin{lemma}\label{scalarreduction}
	The operator $D=P_-\mathcal A|_{E_-}$ is unitarily equivalent to
	\[
	T:=\frac12(A+RBR)=I-H-RJR
	\]
	on $\mathcal H$  where $H$ and $J$ were defined in \eqref{HJ} and $A$ and $B$ in \eqref{AB}. More precisely, for every $f\in\mathcal H$,
	\begin{equation}
		D(f,-Rf)=(Tf,-RTf).
		\label{eq:mit-cr-scalar-reduction}
	\end{equation}
\end{lemma}

\begin{proof}
	Recall that
	\[
	P_-(u,v)=\frac12(u-Rv,v-Ru).
	\]
	Since $\mathcal A(f,-Rf)=(Af,-BRf)$, we obtain
	\[
	D(f,-Rf) { = P_- \mathcal{A}(f,-Rf) }
	=\frac12(Af+RBRf,-BRf-RAf)
	=(Tf,-RTf),
	\]
	where we used $R^2=I$. Thus $DU_-=U_-T$, with
	$U_-f=2^{-1/2}(f,-Rf)$. Since $U_-$ is unitary from
	$\mathcal H$ onto $E_-$, the result follows.
\end{proof}

\medskip\noindent
Hence $D$ and $T$ are unitarily equivalent.
Therefore, $\widehat L$ is Fredholm if and only if $T$ is Fredholm,
and in that case
\[
\operatorname{ind}\widehat L=\operatorname{ind}D
=\operatorname{ind}T.
\]

\subsection{A Fredholm property}

By definition of $T$ and \eqref{HJ} we have
\[
(Tf)(x)
=
f(x)-\int_x^1\frac{f(t)}{t}\,dt
-(1-x)\int_0^x\frac{f(t)}{(1-t)^2}\,dt\,.
\]
Define the bounded operator $T_0$ on $\mathcal H=L^2(0,1)$ by
\[
(T_0f)(x)=
\begin{cases}
	\displaystyle
	f(x)-\int_x^{1/2}\frac{f(t)}{t}\,dt\,,
	&0<x<\frac12\,,\\[1em]
	\displaystyle
	f(x)-(1-x)\int_{1/2}^x\frac{f(t)}{(1-t)^2}\,dt\,,
	&\frac12<x<1\,.
\end{cases}
\]
The operator $T-T_0$ has a bounded integral kernel,
and is therefore Hilbert--Schmidt, hence compact.

\begin{lemma}
	The operator $T_0$ is injective.
\end{lemma}

\begin{proof}
	Let $f$ such that $T_0f=0\,$. On $(0,1/2)\,$, differentiation in the sense of distributions gives
	$f'=-f/x\,$, hence $(xf)'=0\,$. Thus $f(x)=c/x$ almost everywhere;
	since $f\in L^2(0,1/2)\,$, we have $c=0\,$.
	On $(1/2,1)\,$, differentiation combined with the integral equation $T_0(f)(x) = 0$ gives $f'=0\,$, while again the integral
	equation $T_0(f)(x) = 0$ gives $f(\frac{1}{2}^+)=0 \,$. Thus $f=0\,$.
\end{proof}

\begin{lemma}
	We have
	\[
	\operatorname{Ran}T_0
	=
	\left\{g\in\mathcal H:
	\int_0^{1/2}g(x)\,dx=0\right\}\,.
	\]
\end{lemma}

\begin{proof}
	By Fubini's theorem,
	\[
	\int_0^{1/2}(T_0f)(x)\,dx
	=
	\int_0^{1/2}f(x)\,dx
	-
	\int_0^{1/2}\frac{f(t)}{t}
	\left(\int_0^t dx\right)dt
	=0\,.
	\]
	Conversely, let $g\in\mathcal H$ satisfy
	$\int_0^{1/2}g=0\,$.
	As in \cite[Lemma~7.1]{GGHN2026}, we seek an explicit inverse
	on the range. Solving the associated first-order equations
	suggests
	\[
	f(x)=
	\begin{cases}
		\displaystyle
		g(x)-\frac1x\int_0^x g(t)\,dt\,,
		&0<x<\frac12\,,\\[1em]
		\displaystyle
		g(x)+\int_{1/2}^x\frac{g(t)}{1-t}\,dt\,,
		&\frac12<x<1\,.
	\end{cases}
	\]
	Hardy's inequality and its reflected adjoint form imply
	$f\in\mathcal H$. Now, we check that $T_0f=g\,$.
	
	\par\noindent
	On $(0,1/2)\,$, set $F(x)=\int_0^x g(t)\,dt\,$.
	Since $F(1/2)=0$ and
	\[
	\left(\frac{F(x)}x\right)'=\frac{f(x)}x\,,
	\]
	we obtain
	\[
	\int_x^{1/2}\frac{f(t)}t\,dt=-\frac{F(x)}x\,,
	\]
	and hence $(T_0f)(x)=g(x)\,$.

	\par\noindent
	On $(1/2,1)\,$, set
	$Q(x)=\int_{1/2}^x g(t)/(1-t)\,dt\,$.
	Since $Q(1/2)=0$ and
	\[
	\left(\frac{Q(x)}{1-x}\right)'
	=\frac{f(x)}{(1-x)^2}\,,
	\]
	we obtain
	\[
	(1-x)\int_{1/2}^x\frac{f(t)}{(1-t)^2}\,dt=Q(x)\,,
	\]
	which again gives $(T_0f)(x)=g(x)$\,.
\end{proof}

\medskip\noindent
Thus $T_0$ is injective, and its range is the kernel of the nonzero
continuous linear functional $g\mapsto\int_0^{1/2}g(x)\,dx$,
hence is closed and of codimension one. Since $T-T_0$ is compact, $T$ is Fredholm and
\[
\operatorname{ind}T=\operatorname{ind}T_0=-1\,.
\]


\subsection{Compact perturbation and Fredholm index}

We conclude by a compact-perturbation argument as in
\cite[Section~7]{GGHN2026}.
Set
\[
\mu_n=
\begin{cases}
	\pi n+3\pi/4\,,& n\ge0\,,\\
	\pi n-\pi/4\,,& n\le-1\,,
\end{cases}
\]
and define $L^{(0)}$ by replacing $\lambda_{1,n}$ with $\mu_n$
in the second component of $L$ (see \eqref{eq:mit-cr-L-Phi}), leaving the first unchanged.
By \eqref{eq:AsymptAKNS},
\begin{equation}\label{eq:deltan}
(\delta_n)_{n\in\mathbb Z}
:=(\lambda_{1,n}-\mu_n)_{n\in\mathbb Z}
\in\ell^2(\mathbb Z)\,.
\end{equation}
As explained above, the trigonometric vector associated with $\mu_n$
is $e_n^-$ for $n\ge0$ and $e_{n-1}^-$ for $n\le-1$.
Thus $L^{(0)}$ is obtained from $\widehat L$ by deleting
the coordinate $\langle\mathcal A v,e_{-1}^-\rangle_X$
and reindexing the remaining coordinates.
The coordinate projection is surjective with a one-dimensional
kernel, so it has index $1$.
Since $\operatorname{ind}\widehat L=-1$, it follows that
$L^{(0)}$ is Fredholm with
\[
\operatorname{ind}L^{(0)}
=\operatorname{ind}\widehat L+1=0\,.
\]
On the other hand, as a consequence of the mean value theorem, the trigonometric vectors associated with $\lambda_{1,n}$ and
$\mu_n$ differ in $X$-norm by at most $2|\delta_n|$. 
Since $\mathcal A$ is bounded, we obtain from \eqref{eq:deltan}
\[
\|L-L^{(0)}\|_{\mathrm{HS}}^2
\le 4\|\mathcal A\|^2
\sum_{n\in\mathbb Z}|\delta_n|^2<\infty\,.
\]
Hence $L-L^{(0)}$ is compact, and $L$ is Fredholm of index zero.
Since $S=\mathcal D L$ and $\mathcal D$ is an isomorphism,
$S$ is also Fredholm of index zero.
We therefore obtain the following result.

\begin{lemma}[Fredholm property for the pair $(0,1)$]
	\label{prop:mit-cr-closed-range}
	For the MIT bag model,
	the Fr\'echet differential
	\[
	D_{(p,q)}\mathcal S_{0,1}(0,0):
	L^2((0,1);\mathbb R^2)\longrightarrow
	\ell^2(\mathbb Z;\mathbb R)\times\ell^2(\mathbb Z;\mathbb R)\,,
	\]
	is Fredholm of index zero. 
\end{lemma}

\medskip\noindent
Together with the injectivity established in the preceding section,
this achieves the proof of Proposition \ref{iso}.
Since $\mathcal S_{0,1}$ is real-analytic by
\cite[Theorem~3.3]{Serier2006}, the analytic inverse function theorem
implies that it is a real-analytic diffeomorphism between
a neighborhood of zero in $X$ and a neighborhood of
$\mathcal S_{0,1}(0,0)$ in $Y$. This completes the proof of Theorem \ref{thm:main-MIT-01}.

\begin{rem}
In the three-dimensional MIT bag model on the ball, the two lowest
positive values of the effective angular momentum are $\kappa=1$
and $\kappa=2$. It is therefore natural to ask whether the spectra
of $H_1(p,q)$ and $H_2(p,q)$ uniquely determine $(p,q)$ in a
neighborhood of zero in $L^2((0,1);\mathbb R^2)$.
For this pair, both radial operators are singular at the origin.
Extending our approach to this case would involve lengthy symbolic
calculations requiring a computer algebra system.
\end{rem}




\noindent
D.~Gobin\\
Laboratoire de Mathématiques Jean Leray, UMR CNRS 6629\\
Nantes Université\\
F-44000 Nantes, France\\
\emph{e-mail}: damien.gobin@univ-nantes.fr

\medskip

\noindent
B.~Grébert\\
Laboratoire de Mathématiques Jean Leray, UMR CNRS 6629\\
Nantes Université\\
F-44000 Nantes, France\\
\emph{e-mail}: benoit.grebert@univ-nantes.fr

\medskip

\noindent
B.~Helffer\\
Laboratoire de Mathématiques Jean Leray, UMR CNRS 6629\\
Nantes Université\\
F-44000 Nantes, France\\
\emph{e-mail}: Bernard.Helffer@univ-nantes.fr

\medskip

\noindent
F.~Nicoleau\\
Laboratoire de Mathématiques Jean Leray, UMR CNRS 6629\\
Nantes Université\\
F-44000 Nantes, France\\
\emph{e-mail}: francois.nicoleau@univ-nantes.fr
	
\end{document}